\documentclass[11pt]{article}
\usepackage{amsmath,cite,tikz,multirow,amssymb,mathrsfs,latexsym,bm,inputenc}
\begin{document}
\newcommand{\qed}{\hphantom{.}\hfill $\Box$\medbreak}
\newcommand{\proof}{\noindent{\bf Proof \ }}
\newtheorem{Theorem}{Theorem}[section]
\newtheorem{Lemma}[Theorem]{Lemma}
\newtheorem{Corollary}[Theorem]{Corollary}
\newtheorem{Remark}[Theorem]{Remark}
\newtheorem{Example}[Theorem]{Example}
\newtheorem{Definition}[Theorem]{Definition}
\newtheorem{Construction}[Theorem]{Construction}

\thispagestyle{empty}
\renewcommand{\thefootnote}{\fnsymbol{footnote}}

\title{\bf Group divisible $(K_4-e)$-packings with any minimum leave
\footnote{Supported by NSFC under Grant $11431003$, $11471032$, and Fundamental Research Funds for the Central Universities under Grant $2016$JBM$071$, $2016$JBZ$012$}}
\author{Yufeng Gao$^{1,2}$, Yanxun Chang$^{1}$, Tao Feng$^{1}$\\
$^{1}$Institute of Mathematics, Beijing Jiaotong University, Beijing 100044, P. R. China\\
$^{2}$College of Mathematics, Tonghua Normal University, Tonghua 134000, P. R. China\\
yxchang@bjtu.edu.cn, tfeng@bjtu.edu.cn}

\date{ }
\maketitle

\noindent {\bf Abstract:} A decomposition of $K_{n(g)}\setminus L$, the complete $n$-partite equipartite graph with a subgraph $L$ (called the leave) removed, into edge disjoint copies of a graph $G$ is called a maximum group divisible packing of $K_{n(g)}$ with $G$ if $L$ contains as few edges as possible. We examine all possible minimum leaves for maximum group divisible $(K_4-e)$-packings. Necessary and sufficient conditions are established for their existences.




\section{Introduction}

Let $H$ be a simple graph and $\cal G$ be a set of simple graphs. An {\em $(H, {\cal G})$-packing} is a pair $(X, {\cal B})$ where $X$ is the
vertex set of $H$ and ${\cal B}$ is a collection of subgraphs (called {\em blocks}) of $H$, such that each block is isomorphic to a graph of $\cal G$, and each edge of $H$ is contained in at most one block of $\cal B$. If $\cal G$ contains a single graph $G$, we speak of an $(H, G)$-packing.

When $H$ is the complete multipartite graph $K_{u_1,u_2,\ldots,u_t}$, a $( K_{u_1,u_2,\ldots,u_t},G)$-packing is said to be a {\em group divisible packing}, written as a $G$-GDP of type $\{u_1,u_2,\ldots,u_t\}$. Each partite set is called a {\em group}, and the multiset $\{u_1,u_2,\ldots,u_t\}$ is called the {\em group type} (or {\em type}) of the GDP. For simplicity, we use $g_1^{n_1}g_2^{n_2}\cdots g_r^{n_r}$ to denote $n_i$ occurrences of $g_i$ for $1\leq i\leq r$ in group-types.

An $(H, G)$-packing $(X, {\cal B})$ is called {\em maximum} if there does not exist any $(H, G)$-packing $(X, {\cal A})$ with $|{\cal B}|<|{\cal A}|$. The {\em leave graph} (or {\em leave}) of an $(H, G)$-packing $(X,{\cal B})$ is the graph $(X,L)$ consisting of those edges of $H$ which are included in none of blocks of $\cal B$. If the leave of an $(H, G)$-packing is null, then such a packing is maximum, and referred to as an {\em $(H, G)$-design}. A $(K_{u_1,u_2,\ldots,u_t}, G)$-design is often said to be a {\em group divisible design}, denoted by a $G$-GDD of type $\{u_1,u_2,\ldots,u_t\}$.

Consider the necessary conditions for the existence of a $G$-GDP of type $g^n$ with leave $(X,L)$: (1) the size of $X$ is $gn$; (2) the number of edges of $L$ must be congruent to $\binom{n}{2}g^2\ ({\rm mod}\ |E(G)|)$; (3) by counting blocks involving a specified vertex, we have every vertex in leave $(X,L)$ has degree congruent to $g(n-1)$ modulo $d$, where $d$ is the greatest common divisor of the degrees of the vertices in $G$.

We always write $G$-MGDP instead of maximum $G$-GDP. Billington and Lindner \cite{bl} examined all possible leaves among $K_3$-MGDPs of type $g^n$, and showed that

\begin{Theorem}\label{thm:K3}\emph{\cite{bl}}
Let $n\geq 3$. There exists a $K_3$-MGDP of type $g^n$ with leave $(X,L)$ if and only if $(X,L)$ is a simple graph over $gn$ vertices with the least number of edges satisfying that $n(n-1)g^2-2|L| \equiv 0 \pmod{3}$ and every vertex has degree congruent to $g(n-1) \pmod{2}$.
\end{Theorem}

Hu, Chang and Feng \cite{hcf} considered all possible leaves among $(K_3+e)$-MGDPs of type $g^n$, where $K_3+e$ is a triangle with a pendent edge.

\begin{Theorem}\label{thm:K3}\emph{\cite{hcf}}
Let $n\geq 3$. There exists a $(K_3+e)$-MGDP of type $g^n$ with leave $(X,L)$ if and only if $(X,L)$ is a simple graph over $gn$ vertices with the least number of edges satisfying that $n(n-1)g^2-2|L| \equiv 0 \pmod{8}$, except for $(g,n)=(1,4)$. A $(K_3+e)$-MGDP of type $1^4$ admits a path of length $2$ as the unique leave.
\end{Theorem}

This paper concerns about the existence of a $(K_4-e)$-MGDP of type $g^n$ with any possible minimum leave, where $K_4-e$ is the graph with vertices $a$, $b$, $c$, $d$ and edges $ab$, $ac$, $ad$, $bc$, $bd$, denoted by $[a,b,c-d]$. Denote by $D(g^n)$ the number of blocks in a $(K_4-e)$-MGDP of type $g^n$. Obviously $D(g^n)\leq \lfloor n(n-1)g^2/10\rfloor$. If $D(g^n)$ could attain $\lfloor n(n-1)g^2/10\rfloor$, then the number of edges in the leave would be
\begin{eqnarray*}
n(n-1)g^2/2-5\lfloor n(n-1)g^2/10\rfloor&=&
\left\{
\begin{array}{lll}
0, &  n(n-1)g^2\equiv 0\ ({\rm mod}\ 10),\\
1, &  n(n-1)g^2\equiv 2\ ({\rm mod}\ 10),\\
2, &  n(n-1)g^2\equiv 4\ ({\rm mod}\ 10),\\
3, &  n(n-1)g^2\equiv 6\ ({\rm mod}\ 10),\\
4, &  n(n-1)g^2\equiv 8\ ({\rm mod}\ 10),\\
\end{array}
\right.
\end{eqnarray*}
which can be used to yield all possible minimal leaves listed in Table \ref{possible minimal leaves}, where $E_i$, $i=1,2,3,4$, is any graph with $i$ edges. Specifically, $E_{2}$ could be one of the following two graphs:

\begin{table}
\centering
\begin{tabular}{cc|ccccc}\\\hline
  & & & $n\ ({\rm mod}\ 5)$  \\
  ~~~ & ~~~~~~~~ & ~~~~0~~~~ & ~~~~1~~~~ & ~~~~2~~~~ & ~~~~3~~~~ & ~~~~4~~~~ \\
  \hline
  ~ & 0 & $\emptyset$ & $\emptyset$ & $\emptyset$ & $\emptyset$ & $\emptyset$  \\
  ~ & 1 & $\emptyset$ & $\emptyset$ & $E_{1}$ & $E_{3}$ & $E_{1}$ \\
  $g\ ({\rm mod}\ 5)$ & 2 & $\emptyset$ & $\emptyset$ & $E_{4}$ & $E_{2}$ & $E_{4}$ \\
  ~ & 3 & $\emptyset$ & $\emptyset$ & $E_{4}$ & $E_{2}$ & $E_{4}$ \\
  ~ & 4 & $\emptyset$ & $\emptyset$ & $E_{1}$ & $E_{3}$ & $E_{1}$ \\
  \hline
\end{tabular}\caption{possible leaves among $(K_{4}-e)$-MGDPs of type $g^n$}\label{possible minimal leaves}
\end{table}

\begin{center}
\begin{tabular}{cc}\vspace{0.5cm}
$E_{2,1}$:~~
  \begin{tikzpicture}
    \draw (0,0)--(0.5,0)--(1,0);
    \fill (0,0)circle (1pt);
    \fill (0.5,0)circle (1pt);
    \fill (1,0)circle (1pt);
    \end{tikzpicture}~~&
$E_{2,2}$:~~
  \begin{tikzpicture}
    \draw (0,0)--(0.5,0) (1,0)--(1.5,0);
    \fill (0,0)circle (1pt);
    \fill (0.5,0)circle (1pt);
    \fill (1,0)circle (1pt);
    \fill (1.5,0)circle (1pt);
    \end{tikzpicture}
\end{tabular}
\end{center}
$E_{3}$ could be one of the following five graphs:
\begin{center}
\begin{tabular}{ccc}\vspace{0.5cm}
$E_{3,1}$:~~
\begin{tikzpicture}
    \draw (0,0)--(0.5,0) (0,0)--(0.25,0.4) (0.5,0)--(0.25,0.4);
    \fill (0,0)circle (1pt);
    \fill (0.5,0)circle (1pt);
    \fill (0.25,0.4)circle (1pt);
\end{tikzpicture}~~&
$E_{3,2}$:~~
  \begin{tikzpicture}
    \draw (0,0.5)--(0.25,0.3)--(0.5,0.5) (0.25,0.3)--(0.25,0);
    \fill (0,0.5)circle (1pt);
    \fill (0.25,0.3)circle (1pt);
    \fill (0.5,0.5)circle (1pt);
    \fill (0.25,0)circle (1pt);
    \end{tikzpicture}~~&
$E_{3,3}$:~~
  \begin{tikzpicture}
    \draw (0,0.15)--(0.5,0.15)--(1,0.15)--(1.5,0.15);
    \fill (0,0.15)circle (1pt);
    \fill (0.5,0.15)circle (1pt);
    \fill (1,0.15)circle (1pt);
    \fill (1.5,0.15)circle (1pt);
    \end{tikzpicture}~~\\ \vspace{0.5cm}
$E_{3,4}$:~~
  \begin{tikzpicture}
    \draw (0,0.15)--(0.5,0.15) (1,0.15)--(1.5,0.15) (2,0.15)--(1.5,0.15);
    \fill (0,0.15)circle (1pt);
    \fill (0.5,0.15)circle (1pt);
    \fill (1,0.15)circle (1pt);
    \fill (1.5,0.15)circle (1pt);
    \fill (2,0.15)circle (1pt);
    \end{tikzpicture}~~&
$E_{3,5}$:~~
  \begin{tikzpicture}
    \draw (0,0.15)--(0.5,0.15) (1,0.15)--(1.5,0.15) (2,0.15)--(2.5,0.15);
    \fill (0,0.15)circle (1pt);
    \fill (0.5,0.15)circle (1pt);
    \fill (1,0.15)circle (1pt);
    \fill (1.5,0.15)circle (1pt);
    \fill (2,0.15)circle (1pt);
    \fill (2.5,0.15)circle (1pt);
\end{tikzpicture}~~
\end{tabular}
\end{center}
$E_{4}$ could be one of the following eleven graphs:
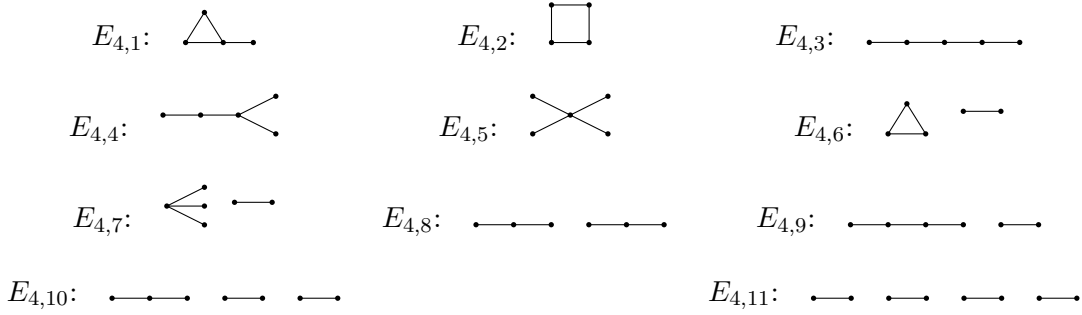
\begin{center}
\begin{tabular}{ccc}\vspace{0.5cm}
$E_{4,1}$:~~
\begin{tikzpicture}
    \draw (0,0)--(0.5,0) (0,0)--(0.25,0.4) (0.5,0)--(0.25,0.4) (0.5,0)--(0.9,0);
    \fill (0,0)circle (1pt);
    \fill (0.5,0)circle (1pt);
    \fill (0.25,0.4)circle (1pt);
    \fill (0.9,0)circle (1pt);
\end{tikzpicture}~~&
$E_{4,2}$:~~
\begin{tikzpicture}
    \draw (0,0)--(0.5,0)--(0.5,0.5)--(0,0.5)--cycle;
    \fill (0,0)circle (1pt);
    \fill (0.5,0)circle (1pt);
    \fill (0.5,0.5)circle (1pt);
    \fill (0,0.5)circle (1pt);
\end{tikzpicture}~~&
$E_{4,3}$:~~
  \begin{tikzpicture}
    \draw (0,0.15)--(0.5,0.15)--(1,0.15)--(1.5,0.15)--(2,0.15);
    \fill (0,0.15)circle (1pt);
    \fill (0.5,0.15)circle (1pt);
    \fill (1,0.15)circle (1pt);
    \fill (1.5,0.15)circle (1pt);
    \fill (2,0.15)circle (1pt);
    \end{tikzpicture}~~\\ \vspace{0.5cm}
$E_{4,4}$:~~
    \begin{tikzpicture}
    \draw  (0,0.25)--(0.5,0.25) (0.5,0.25)--(1,0.25)
    (1,0.25)--(1.5,0.5) (1,0.25)--(1.5,0);
    \fill (0,0.25)circle (1pt);
    \fill (0.5,0.25)circle (1pt);
    \fill (1,0.25)circle (1pt);
    \fill (1.5,0.5)circle (1pt);
    \fill (1.5,0)circle (1pt);
\end{tikzpicture}~~&
$E_{4,5}$:~~
  \begin{tikzpicture}
    \draw (0,0.5)--(0.5,0.25) (0,0)--(0.5,0.25) (1,0.5)--(0.5,0.25)
    (1,0)--(0.5,0.25);
    \fill (0,0.5)circle (1pt);
    \fill (0.5,0.25)circle (1pt);
    \fill (0,0)circle (1pt);
    \fill (1,0.5)circle (1pt);
    \fill (1,0)circle (1pt);
   \end{tikzpicture}~~&
$E_{4,6}$:~~
\begin{tikzpicture}
    \draw (0,0)--(0.5,0) (0,0)--(0.25,0.4) (0.5,0)--(0.25,0.4) (1,0.3)--(1.5,0.3);
    \fill (0,0)circle (1pt);
    \fill (0.5,0)circle (1pt);
    \fill (0.25,0.4)circle (1pt);
    \fill (1,0.3)circle (1pt);
    \fill (1.5,0.3)circle (1pt);
\end{tikzpicture}~~\\ \vspace{0.5cm}
%
$E_{4,7}$:~~

\begin{tikzpicture}
    \draw (0,0.25)--(0.5,0.5) (0,0.25)--(0.5,0.25) (0,0.25)--(0.5,0) (0.9,0.3)--(1.4,0.3);
    \fill (0,0.25)circle (1pt);
    \fill (0.5,0.5)circle (1pt);
    \fill (0.5,0.25)circle (1pt);
    \fill (0.5,0)circle (1pt);
        \fill (0.9,0.3)circle (1pt);
    \fill (1.4,0.3)circle (1pt);
\end{tikzpicture}~~&
$E_{4,8}$:~~
  \begin{tikzpicture}
    \draw (0,0.15)--(0.5,0.15)--(1,0.15) (1.5,0.15)--(2,0.15) (2,0.15)--(2.5,0.15);
    \fill (0,0.15)circle (1pt);
    \fill (0.5,0.15)circle (1pt);
    \fill (1,0.15)circle (1pt);
    \fill (1.5,0.15)circle (1pt);
    \fill (2,0.15)circle (1pt);
    \fill (2.5,0.15)circle (1pt);
    \end{tikzpicture}~~&
$E_{4,9}$:~~
  \begin{tikzpicture}
    \draw (0,0.15)--(0.5,0.15)--(1,0.15)--(1.5,0.15) (2,0.15)--(2.5,0.15);
    \fill (0,0.15)circle (1pt);
    \fill (0.5,0.15)circle (1pt);
    \fill (1,0.15)circle (1pt);
    \fill (1.5,0.15)circle (1pt);
    \fill (2,0.15)circle (1pt);
    \fill (2.5,0.15)circle (1pt);
\end{tikzpicture}~~\\ \vspace{0.5cm}
$E_{4,10}$:~~
  \begin{tikzpicture}
    \draw (0,0.15)--(0.5,0.15)--(1,0.15) (1.5,0.15)--(2,0.15) (2.5,0.15)--(3,0.15);
    \fill (0,0.15)circle (1pt);
    \fill (0.5,0.15)circle (1pt);
    \fill (1,0.15)circle (1pt);
    \fill (1.5,0.15)circle (1pt);
    \fill (2,0.15)circle (1pt);
    \fill (2.5,0.15)circle (1pt);
    \fill (3,0.15)circle (1pt);
\end{tikzpicture}~~&
~&
$E_{4,11}$:~~
  \begin{tikzpicture}
    \draw (0,0.15)--(0.5,0.15) (1,0.15)--(1.5,0.15) (2,0.15)--(2.5,0.15) (3,0.15)--(3.5,0.15);
    \fill (0,0.15)circle (1pt);
    \fill (0.5,0.15)circle (1pt);
    \fill (1,0.15)circle (1pt);
    \fill (1.5,0.15)circle (1pt);
    \fill (2,0.15)circle (1pt);
    \fill (2.5,0.15)circle (1pt);
    \fill (3,0.15)circle (1pt);
    \fill (3.5,0.15)circle (1pt);
\end{tikzpicture}~~
\end{tabular}
\end{center}

Colbourn, Ling and Quattrocchi \cite{clq} established the necessary and sufficient condition for the existence of a $(K_4-e)$-GDD of type $g^n$.

\begin{Theorem}\label{thm:K4-e-GDD}\emph{\cite{clq}}
There exists a $(K_{4}-e)$-GDD of type $g^n$ if and only if $n\geq 3$, $gn(n-1)\equiv 0\pmod{5}$ and $(g,n)\neq(1,5)$.
\end{Theorem}

Hoffman, Lindner, Sharry and Street \cite{hlss} investigated all possible leaves in a $(K_4-e)$-MGDP of type $1^n$. They proved that
\begin{eqnarray*}
D(1^n)&=&
\left\{
\begin{array}{lll}
\lfloor n(n-1)/10\rfloor-1, &  n=5,7,9,\\
\lfloor n(n-1)/10\rfloor, &  {\rm otherwise}.
\end{array}
\right.
\end{eqnarray*}
It follows that for $g=1$ and $n=5,7,9$, all possible leaves among $(K_4-e)$-MGDPs of type $1^n$ do not satisfy Table \ref{possible minimal leaves}. Actually there is only one possible leave for $n=5$; there are 7 possible leaves for $n=7$; there are 36 possible leaves for $n=9$ (see Example 1.1, 1.2 and 1.4 in \cite{hlss}).

\begin{Theorem} \label{thm:1^n}\emph{\cite{hlss}}
Suppose that $L$ is any simple graph on $n$ vertices with $\binom{n}{2}-5D(1^n)$ edges. There exists a $(K_4-e)$-MGDP of type $1^n$ with leave $L$ for any $n\geq 4$ with the exception of $n\in\{5,7,8,9\}$. A $(K_4-e)$-MGDP of type $1^8$ can only admit $E_{3,j}$ as its leave for $j=3,4,5$.
\end{Theorem}

This paper continues the study in \cite{hlss} to figure out all possible leaves in $(K_4-e)$-MGDPs of type $g^n$ for any $g>1$. Throughout this paper we always assume that $I_u=\{0,1,\ldots,u-1\}$.

\begin{Example} \label{ex:2^3}
We here construct a $(K_{4}-e)$-MGDP of type $2^3$ on $I_6$ with groups $\{i,i+3\}$, $i=0,1,2$, which admits two disjoint edges $\{0,4\}$ and $\{1,3\}$ as its leave. Its two blocks are $[0,  1,  2-  5]$ and $[3,  4,  2-  5]$.
\end{Example}

As the main result of the paper, we are to prove

\begin{Theorem}\label{thm:main}
Let $n\geq 3$. There exists a $(K_{4}-e)$-MGDP of type $g^n$ with leave $(X,L)$ if and only if $(X,L)$ is a simple graph over $gn$ vertices with the least number of edges satisfying that $n(n-1)g^2-2|L| \equiv 0 \pmod{10}$, except when $(g,n)\in \{(1,5), (1,7), (1,8),(1,9), (2,3),(2,4),(3,3),(4,$ $3)\}$. A $(K_{4}-e)$-MGDP of type $g^3$, $g\in\{2,3\}$, admits two disjoint edges as the unique leave. A $(K_{4}-e)$-MGDP of type $4^3$ admits a triangle as the unique leave. A $(K_{4}-e)$-MGDP of type $2^4$ can only admit $E_{4,j}$ as its leave for $2\leq j\leq 11$.
\end{Theorem}

\begin{Remark}\label{rmk}
It is readily checked that if $(X,L)$ is a simple graph over $gn$ vertices with the least number of edges satisfying that $n(n-1)g^2-2|L| \equiv 0 \pmod{10}$, then $|L|=n(n-1)g^2/2-5\lfloor n(n-1)g^2/10$. Thus Theorem \ref{thm:main} implies $D(g^n)=\lfloor n(n-1)g^2/10\rfloor$.
\end{Remark}

\section{Nonexistence results}

\begin{Lemma} \label{lem:2^3_no}
There is no $(K_{4}-e)$-MGDP of type $2^3$ with leave $E_{2,1}$.
\end{Lemma}

\proof Suppose that there exists a $(K_{4}-e)$-MGDP of type $2^3$ with leave $E_{2,1}$ on $I_6$ with groups $\{i,i+3\}$, $0\leq i\leq 2$. It has $2$ blocks. Denote by ${\cal B}$ the set of its blocks. W.o.l.g, assume that $E_{2,1}$ consists of two edges $\{0,1\}$, $\{1,2\}$, or $\{0,1\}$, $\{1,3\}$. For every $x\in I_6$, denote by $d_{i}(x)$, $i=2,3$, the number of blocks in ${\cal B}$ in which the degree of $x$ is $i$.

If $E_{2,1}$ consists of two edges $\{0,1\}$ and $\{1,2\}$, then
$$2d_2(x)+3d_3(x)=\left\{
\begin{array}{lll}
3, & x\in\{0,2\};\\
2, & x=1;\\
4, & x\in\{3,4,5\}.\\
\end{array}
\right.
$$
Then $d_{2}(x)=0$ and $d_{3}(x)=1$ for $x\in\{0,2\}$; $d_{2}(1)=1$ and $d_3(1)=0$; $d_2(x)=2$ and $d_3(x)=0$ for $x\in\{3,4,5\}$.

If $E_{2,1}$ consists of two edges $\{0,1\}$ and $\{1,3\}$, then
$$2d_2(x)+3d_3(x)=\left\{
\begin{array}{lll}
3, & x\in\{0,3\};\\
2, & x=1;\\
4, & x\in\{2,4,5\}.\\
\end{array}
\right.
$$
Then $d_{2}(x)=0$ and $d_{3}(x)=1$ for $x\in\{0,3\}$; $d_{2}(1)=1$ and $d_3(1)=0$; $d_2(x)=2$ and $d_3(x)=0$ for $x\in\{2,4,5\}$.

Therefore, $\sum_{x=0}^{5}d_{3}(x)=2$. On the other hand, since $|{\cal B}|=2$, $\sum_{x=0}^{5}d_{3}(x)=4$, a contradiction. \qed

\begin{Lemma} \label{lem:3^3_no}
There is no $(K_{4}-e)$-MGDP of type $3^3$ with leave $E_{2,1}$.
\end{Lemma}

\proof Suppose that there exists a $(K_{4}-e)$-MGDP of type $3^3$ with leave $E_{2,1}$ on $I_9$ with groups $\{i,i+3,i+6\}$, $0\leq i\leq 2$. It has $5$ blocks. Denote by ${\cal B}$ the set of its blocks. W.o.l.g., assume that $E_{2,1}$ consists of two edges $\{0,1\}$, $\{1,2\}$, or $\{0,1\}$, $\{1,3\}$. For every $x\in I_9$, denote by $d_{i}(x)$, $i=2,3$, the number of blocks in ${\cal B}$ in which the degree of $x$ is $i$.

If $E_{2,1}$ consists of two edges $\{0,1\}$ and $\{1,2\}$, then
$d_{2}(x)=d_{3}(x)=1$ for $x\in\{0,2\}$; $d_{2}(1)=2$ and $d_3(1)=0$; $d_2(x)=3$ and $d_3(x)=0$, or $d_2(x)=0$ and $d_3(x)=2$ for $x\in\{3,4,5,6,7,8\}$. Analyzing blocks containing 0 and 1, we have $\cal B$ must contain blocks of the form
$$[0,*,*-*],[y_1,y_2,0-1],[y_3,y_4,*-1],$$
or
$$[0,*,*-*],[*,*,0-*],[z_1,z_2,*-1],[z_3,z_4,*-1],[z_5,z_6,*-*].$$
In the former case, $\{y_1,y_2\}$ must be $\{5,8\}$. However, $5$ and $8$ are in the same group, a contradiction. In the latter case, w.o.l.g., we can assume that $\{z_1,z_2\}=\{3,5\}$ and $\{z_3,z_4\}=\{6,8\}$. Thus $d_3(x)=2$ for $x\in\{3,5,6,8\}$, which yields $\{z_5,z_6\}=\{3,6\}$, a contradiction.

If $E_{2,1}$ consists of two edges $\{0,1\}$ and $\{1,3\}$, then $2d_2(1)+3d_3(1)=4$, which yields $d_{2}(1)=2$ and $d_3(1)=0$. Thus $\cal B$ contains blocks of the form
$[y_1,y_2,1-*]$ and $[y_3,y_4,1-*]$, where $\{y_1,y_2,y_3,y_4\}=\{2,5,6,8\}$. However, $2,5$ and $8$ are in the same group, a contradiction. \qed

\begin{Lemma} \label{lem:2^4_no}
There is no $(K_{4}-e)$-MGDP of type $2^4$ with leave $E_{4,1}$.
\end{Lemma}

\proof Suppose that there exists a $(K_{4}-e)$-MGDP of type $2^4$ with leave $E_{4,1}$ on $I_8$ with groups $\{i,i+4\}$, $0\leq i\leq 3$. It has $4$ blocks. Denote by ${\cal B}$ the set of its blocks. W.o.l.g, assume that $E_{4,1}$ consists of four edges $\{0,1\}$, $\{1,2\}$, $\{0,2\}$, $\{0,5\}$, or $\{0,1\}$, $\{1,2\}$, $\{0,2\}$, $\{0,3\}$. For every $x\in I_8$, denote by $d_{i}(x)$, $i=2,3$, the number of blocks in ${\cal B}$ in which the degree of $x$ is $i$.

If $E_{4,1}$ consists of four edges $\{0,1\}$, $\{1,2\}$, $\{0,2\}$ and $\{0,5\}$, then
$d_{2}(0)=0$ and $d_3(0)=1$; $d_{2}(x)=2$ and $d_{3}(x)=0$ for $x\in\{1,2\}$; $d_{2}(5)=d_{3}(5)=1$; $d_2(x)=3$ and $d_3(x)=0$, or $d_2(x)=0$ and $d_3(x)=2$ for $x\in\{3,4,6,7\}$. Since the unique block containing 0 is $[0,6,3-7]$, we have $d_2(x)=3$ and $d_3(x)=0$ for $x\in\{3,7\}$. Thus $\sum_{x=0}^{7}d_{2}(x)\geq 11$. On the other hand, since $|{\cal B}|=4$, $\sum_{x=0}^{7}d_{2}(x)=8$, a contradiction.

If $E_{4,1}$ consists of four edges $\{0,1\}$, $\{1,2\}$, $\{0,2\}$ and $\{0,3\}$, then
$d_{2}(0)=0$ and $d_3(0)=1$; $d_{2}(x)=2$ and $d_{3}(x)=0$ for $x\in\{1,2\}$; $d_{2}(3)=d_{3}(3)=1$; $d_2(x)=3$ and $d_3(x)=0$, or $d_2(x)=0$ and $d_3(x)=2$ for $x\in\{4,5,6,7\}$. Analyzing blocks containing 0 and 3, we have $\cal B$ consists of blocks of the form $$[0,*,*-*],[3,y_1,y_2-y_3],[y_4,y_5,*-3],[*,*,*-*].$$
Since $1$ and $2$ cannot be vertices with degree 3 in any block, examining edges $\{1,3\}$ and $\{2,3\}$, we have $\{y_2,y_3\}=\{1,2\}$, which yields $\{y_1,y_4,y_5\}=\{4,5,6\}$. So $d_2(x)=0$ and $d_3(x)=2$ for $x\in\{4,5,6\}$. On the other hand, there is at least one vertex taken from the set $\{5,6\}$ occurring in the block of the form $[0,*,*-*]$ as a vertex with degree 2, a contradiction. \qed

\begin{Lemma} \label{lem:4^3_no}
There is no $(K_{4}-e)$-MGDP of type $4^3$ with leave $E_{3,j}$ for any $j\in\{2,3,4,5\}$.
\end{Lemma}

\proof Suppose that there exists a $(K_{4}-e)$-MGDP of type $4^3$ with leave $L$ on $I_{12}$ with groups $G_i=\{i,i+3,i+6,i+9\}$, $0\leq i\leq 2$, where $L=E_{3,j}$ for some $j\in\{2,3,4,5\}$. It has $9$ blocks. Denote by ${\cal B}$ the set of its blocks. For $x\in I_{12}$, denote by $d_{i}(x)$, $i=2,3$, the number of blocks in ${\cal B}$ in which the degree of $x$ is $i$.

If $x$ is an isolated vertex in $L$, then $d_{2}(x)=1$ and $d_{3}(x)=2$ (such element $x$ is called an {\em a-element}), or $d_2(x)=4$ and $d_3(x)=0$ (such element $x$ is called a {\em b-element}). Denote by $\alpha$, $\beta$ the number of a-elements and b-elements in $L$, respectively.

If $x$ is a vertex with degree one in $L$, then $d_{2}(x)=2$ and $d_{3}(x)=1$. Denote by $\theta$ the number of such kind of vertices in $L$.

If $x$ is a vertex with degree two in $L$, then $d_{2}(x)=3$ and $d_{3}(x)=0$ (such element $x$ is called a {\em c-element}), or $d_2(x)=0$ and $d_3(x)=2$ (such element $x$ is called a {\em d-element}). Denote by $\gamma$, $\delta$ the number of c-elements and d-elements in $L$, respectively.

If $x$ is a vertex with degree three in $L$, then $d_{2}(x)=1$ and $d_{3}(x)=1$. Denote by $\rho$ the number of such kind of vertices in $L$.

There is no vertex with degree greater than three in $L$, so $\alpha+\beta+\theta+\gamma+\delta+\rho=12$. Since each block of
${\cal B}$ contains exactly two vertices with degree $3$, we have
$2\alpha+\theta+2\delta+\rho=2|{\cal B}|=18$. Solving these two equations according to the type of $L$, we have the following table.

\begin{center}
\begin{tabular}{c|llllll}
 & $\alpha$ & $\beta$ & $\theta$ & $\gamma$ & $\delta$ & $\rho$ \\\hline
$E_{3,2}$ & $7$ & $1$ & $3$ & $0$ & $0$ & $1$ \\\hline
 & $8$ & $0$ & $2$ & $2$ & $0$ & $0$ \\
$E_{3,3}$  & $7$ & $1$ & $2$ & $1$ & $1$ & $0$ \\
  & $6$ & $2$ & $2$ & $0$ & $2$ & $0$ \\\hline
$E_{3,4}$  & $7$ & $0$ & $4$ & $1$ & $0$ & $0$ \\
  & $6$ & $1$ & $4$ & $0$ & $1$ & $0$ \\\hline
$E_{3,5}$  & $6$ & $0$ & $6$ & $0$ & $0$ & $0$ \\
\end{tabular}
\end{center}

(1) If all edges in $L$ are from the same two groups, then $L$ contains four isolated vertices coming from the same group. W.o.l.g., the four isolated vertices are from $G_2$. Since $\beta\leq 2$, there are at most two b-elements in $G_2$.

(1.1) If the four vertices in $G_2$ are all a-elements, then $d_3(x)=2$ for any $x\in G_2$. It follows that $\cal B$ consists of blocks of the form \begin{center}
\begin{tabular}{llllll}
$[2,*,*-*]$,&$[2,*,*-*]$,&$[5,*,*-*]$,&$[5,*,*-*]$,&$[8,*,*-*]$,\\
$[8,*,*-*]$,&$[11,*,*-*]$,&$[11,*,*-*]$,& $[*,*,*-*]$.
\end{tabular}
\end{center}
A contradiction occurs since $d_2(x)=1$ for any $x\in G_2$.

(1.2) If there are three a-elements and one b-element in $G_2$, then w.l.o.g., $2,5$ and $8$ are a-elements. It follows that $\cal B$ consists of blocks of the form \begin{center}
\begin{tabular}{llllll}
$[2,*,*-*]$,&$[2,*,*-*]$,&$[5,*,*-*]$,&$[5,*,*-*]$,&$[8,*,*-*]$,\\
$[8,*,*-*]$,&$[*,*,*-*]$,&$[*,*,*-*]$,& $[*,*,*-*]$.
\end{tabular}
\end{center}
A contradiction occurs since $11$ is a b-element.

(1.3) If there are two a-elements and two b-elements in $G_2$, then w.l.o.g., $2,5$ are a-elements and $8,11$ are b-elements. It follows that w.l.o.g., $\cal B$ consists of blocks of the form \begin{center}
\begin{tabular}{llllll}
$[2,*,*-*]$,&$[2,*,*-*]$,&$[5,*,*-*]$,&$[5,*,*-*]$,&$[*,*,8-11]$,\\
$[*,*,8-11]$,&$[*,*,8-11]$,&$[*,*,8-x_1]$,& $[x_4,x_5,x_2-x_3]$,
\end{tabular}
\end{center}
where $(x_1,x_2,x_3)=(11,2,5)$ or $(2,5,11)$. Note that $\beta=2$ only if $L=E_{3,3}$ and $\delta=2$. Thus w.o.l.g, assume that $L$ consists of edges $\{0,1\}$, $\{1,3\}$ and $\{3,4\}$. Then $1$ and $3$ are d-elements, i.e., $d_3(1)=d_3(3)=2$. Since $\{1,2\}$, $\{1,5\}$, $\{1,8\}$ and $\{1,11\}$ must be covered by blocks of $\cal B$, we have $x_4=1$. Since $\{3,2\}$, $\{3,5\}$, $\{3,8\}$ and $\{3,11\}$ must be covered by blocks of $\cal B$, we have $x_5=3$. Thus $\{1,3\}$ is covered by the last block, a contradiction.

(2) If only two edges in $L$ are from the same two groups, then there exists one group in $L$ containing three isolated vertices and one vertex with degree one. W.o.l.g., assume that the group is $G_2$ and the vertex with degree one is $2$. Since $\beta\leq 2$, there are at most two b-elements in $\{5,8,11\}$.

(2.1) If $5,8$ and $11$ are all a-elements, then $d_3(x)=2$ for $x\in\{5,8,11\}$, and $\cal B$ consists of blocks of the form
\begin{center}
\begin{tabular}{llllll}
$[5,*,*-*]$,&$[5,*,*-*]$,&$[8,*,*-*]$,&$[8,*,*-*]$,&$[11,*,*-*]$,\\
$[11,*,*-*]$,&$[2,*,*-*]$,&$[*,*,*-2]$,& $[*,*,*-2]$.
\end{tabular}
\end{center}
A contradiction occurs since $d_2(x)=1$ for $x\in\{5,8,11\}$.

(2.2) If there are two b-elements and one a-element in $\{5,8,11\}$, then $\beta=2$, so $L$ must be $E_{3,3}$. W.l.o.g., assume that $L$ consists of edges $\{0,1\}$, $\{1,3\}$ and $\{3,2\}$. Then $4,7,10$ are a-elements and $1,3$ are d-element. It follows that $\cal B$ consists of blocks of the form \begin{center}
\begin{tabular}{llllll}
$[4,*,*-*]$,&$[4,*,*-*]$,&$[7,*,*-*]$,&$[7,*,*-*]$,&$[10,*,*-*]$,\\
$[10,*,*-*]$,&$[1,*,*-*]$,&$[1,*,*-*]$,& $[*,*,*-*]$.
\end{tabular}
\end{center}
A contradiction occurs since $d_2(x)=1$ for $x\in\{4,7,10\}$.

(2.3) If there are two a-elements and one b-element in $\{5,8,11\}$, then $\beta\geq 1$ and $L$ could be $E_{3,2}$, $E_{3,3}$ or $E_{3,4}$. W.l.o.g., $5,8$ are a-elements and $11$ is a b-element. It follows that $\cal B$ consists of blocks of the form
\begin{center}
\begin{tabular}{llllll}
$[5,*,*-*]$,&$[5,*,*-*]$,&$[8,*,*-*]$,&$[8,*,*-*]$,&$[2,*,*-*]$,\\
$[*,*,2-11]$, & $[*,*,2-11]$,&$[*,*,5-11]$,&$[*,*,8-11]$.
\end{tabular}
\end{center} 

(2.3.1) If $L=E_{3,3}$ and $\delta=2$, or $L=E_{3,4}$ and $\delta=1$, then there exists one d-element $x^*$ in $G_0\cup G_1$ such that $\{x^*,2\}$ is not an edge in $L$. Since $d_2(x^*)=0$, $d_3(x^*)=2$, and $\{x^*,2\}$, $\{x^*,5\}$, $\{x^*,8\}$, $\{x^*,11\}$ must be covered by blocks of $\cal B$, a contradiction.

(2.3.2) If $L=E_{3,3}$ and $\gamma=1$, then there exists one c-element $x^*$ in $G_0\cup G_1$ satisfying $d_2(x^*)=3$ and $d_3(x^*)=0$. Since $11$ is a b-element, $d_2(11)=4$, $d_3(11)=0$, and $\{x^*,11\}$ must be contained in some block of $\cal B$, a contradiction.

(2.3.3) If $L=E_{3,2}$, w.l.o.g., $E_{3,2}$ consists of edges $\{0,1\}$, $\{0,4\}$ and $\{0,2\}$. Note that $11$ is the unique b-element in the current case. Since $3,6$ and $9$ are all a-elements, w.l.o.g., we rewrite $\cal B$ such that $\cal B$ consists of blocks of the form \begin{center}
\begin{tabular}{lllll}
$[3,z_1,*-*]$,&$[3,z_2,z_3-z_4]$,&$[6,z_5,*-*]$,&$[6,*,z_6-z_7]$,&$[9,*,z_8-z_9]$,\\
$[9,z_{10},*-,z_{11}]$,&
$[z_{21},z_{12},3-z_{13}]$,&
$[z_{14},z_{15},6-z_{16}]$,&
$[z_{17},z_{18},z_{19}-z_{20}]$,
\end{tabular}
\end{center}
where $9\in\{z_{13},z_{16},z_{20}\}$. Consider blocks containing $0$. Since $0,3,6$ and $9$ are in the same group, we have $z_{17}=0$, and either $z_{13}=0$ or $z_{16}=0$. W.l.o.g., set $z_{16}=0$, which implies $z_{13}=9$, and then w.o.l.g we can take $z_{14}=5$, $z_{15}=7$, $z_{18}=10$, $z_{19}=8$ and $z_{20}=11$.

Consider blocks containing $8$. If $z_{12}\neq 8$ and $z_{21}\neq 8$, since edges $\{3,8\}$, $\{6,8\}$ and $\{8,9\}$ must be covered by blocks, the number of blocks containing 8 is at least 4. But it is impossible since $8$ is an a-element. Thus w.o.l.g., set $z_{12}=8$, and so $z_5=8$.

Consider blocks containing $5$. Since $5$ is an a-element and $5$ and $8$ are in the same group, w.l.o.g., we can set $z_1=z_8=5$ (this can be done because $3$ and $9$ are symmetric so far).

Consider blocks containing $11$. Since $11$ is the unique b-element and $5,8$ and $11$ are in the same group, we have $z_4=z_7=11$ and $11\in \{z_{9},z_{11}\}$.

Consider blocks containing $2$. Since $2,5,8$ and $11$ are in the same group and edges $\{2,3\}$, $\{2,6\}$ and $\{2,9\}$ must be covered by blocks, we have $z_3=z_6=z_{10}=2$, which yields $z_{11}\neq 11$ and $z_9=11$.

Consider blocks containing $7$. Since $7$ is an a-element and edges $\{0,7\}$, $\{5,7\}$ and $\{6,7\}$ have been covered by some blocks, only $z_2$, $z_{11}$ and $z_{21}$ could be $7$. However, only one of these three values could be $7$, a contradiction.

(3) If any two edges in $L$ are from different two groups, then $L$ cannot be $E_{3,2}$. 
In a $K_4-e$, there is only one edge whose two endpoints both have degree three. We refer to this edge as the {\em distinguished edge} of the $K_4-e$. Since $|{\cal B}|=9$, we have $9$ distinguished edges. If a distinguished edge is from $G_r\cup G_s$, $0\leq r<s\leq 2$, then it is called of type $(r,s)$. Denote by $\eta_1$, $\eta_2$, $\eta_3$ the number of distinguished edges of type $(0,1)$, $(0,2)$ and $(1,2)$, respectively. Then $\eta_1+\eta_2+\eta_3=9$. Counting the number of edges from $G_r\cup G_s$ which are covered by blocks of $\cal B$, $0\leq r<s\leq 2$, we have
$$\left\{
\begin{array}{lll}
(r,s)=(0,1): & \eta_1+2\eta_2+2\eta_3=15;\\
(r,s)=(0,2): & \eta_2+2\eta_1+2\eta_3=15;\\
(r,s)=(1,2): & \eta_3+2\eta_1+2\eta_2=15.\\
\end{array}
\right.
$$
It follows that $\eta_1=\eta_2=\eta_3=3$. This fact yields $\beta=0$. Otherwise, there exists a b-element $x^*$ from some group $G_l$ satisfying $d_2(x^*)=4$. It follows that the four blocks containing $x^*$ contributes four distinguished edges which are all from groups $G_{l'}$ and $G_{l''}$, $l'\neq l$, $l''\neq l$, a contradiction. Therefore, all isolated vertices in $L$ are a-elements. Furthermore, the following two claims are useful.

Claim 1. If $x$ is a vertex with degree one in $L$ and $\{x,y\}$ is an edge in $L$, where $x\in G_r$, $y\in G_s$ and $r\neq s$, then since $d_3(x)=1$, there exists a unique distinguished edge containing $x$, written as $\{x,z\}$. We have that $z$ must come from $G_s$. Otherwise, if $z\in G_t$ and $t\not\in \{r,s\}$, then $\cal B$ contains a block of the form $[x,z,y_1-y_2]$, where $y_1,y_2\in G_{s}$, $y_1\neq y$ and $y_2\neq y$. Note that $d_2(x)=2$. It is readily checked that one cannot get the other two blocks containing $x$, a contradiction.

Claim 2. If $x$ is an isolated vertex in $L$, then $d_3(x)=2$, and so there are two distinguished edges containing $x$, written as $\{x,y\}$ and $\{x,z\}$. We have that $y$ and $z$ must come from different groups. Otherwise, if $y,z\in G_s$ and $x\in G_r$, $r\neq s$, then $\cal B$ contains two blocks of the form $[x,y,u_1-u_2]$ and $[x,z,u_3-u_4]$, where $\{u_1,u_2,u_3,u_4\}=G_t$, $t\not\in\{r,s\}$. Note that $d_2(x)=1$. It is readily checked that one cannot get the last block containing $x$, a contradiction.

W.o.l.g., assume that $E_{3,3}$ consists of edges $\{4,5\}$, $\{5,6\}$ and $\{6,7\}$; $E_{3,4}$ consists of edges $\{0,1\}$, $\{4,5\}$ and $\{5,6\}$; $E_{3,5}$ consists of edges $\{0,1\}$, $\{4,5\}$ and $\{6,8\}$. For each $L=E_{3,j}$, $j=3,4,5$, we have $d_2(4)=2$ and $d_3(4)=1$. It follows that $\cal B$ contains blocks of the form
\begin{center}
\begin{tabular}{lll}
$[4,x_1,x_2-x_3]$,&$[x_4,x_5,x_6-4]$,&$[x_7,x_8,x_9-4]$.
\end{tabular}
\end{center}

(3.1) For $L=E_{3,3}$, due to $\beta=0$, we have $\gamma=2$, and so $5$ and $6$ are c-elements. Since $6$ is a c-element and edge $\{4,6\}$ must be covered be some block, w.o.l.g., set $x_2=6$. Then w.l.o.g., set $x_1=2$, $x_5=8$, $x_8=11$, and $x_3=0$, $x_4=3$, $x_7=9$.

Recall that there are three distinguished edges of type $(0,2)$ in $\cal B$. So far we have gotten two of them, i.e., $\{3,8\}$ and $\{9,11\}$. Since $5$ and $6$ are c-elements, the third one cannot contain $5$ and $6$. Because of Claim 2, the third one cannot be $\{0,8\}$, $\{0,11\}$, $\{3,2\}$, $\{3,11\}$, $\{9,2\}$ and $\{9,8\}$. Thus the third one would be $\{0,2\}$, but $\{0,2\}$ has been covered by the first block containing $4$, a contradiction.

(3.2) For $L=E_{3,4}$, due to $\beta=0$, we have $\gamma=1$, and so $5$ is a c-element. 
Since $4$ is a vertex with degree one in $E_{3,4}$ and $\{4,5\}$ is an edge in $E_{3,4}$, by Claim 1, there exists a unique distinguished edge containing $4$, written as $\{4,z\}$, where $z\in G_2$. Thus w.l.o.g., set $x_1=2$.

Since $0$ is a vertex with degree one in $E_{3,4}$ and $\{0,1\}$ is an edge in $E_{3,4}$, then by Claim 1, there exists a unique distinguished edge containing $0$, written as $\{0,z'\}$, where $z'\in G_1$. Thus $0\not\in\{x_4,x_5,x_7,x_8\}$ (note that $d_3(0)=1$). Since edge $\{0,4\}$ must be covered be some block, w.o.l.g., set $x_2=0$. Then w.l.o.g., set $x_5=8$ and $x_8=11$. It follows that $\{x_3,x_4,x_7\}=\{3,6,9\}$.

If $x_3=6$, then w.o.l.g., set $x_4=3$ and $x_7=9$. Since $6$ is a vertex with degree one in $E_{3,4}$ and $\{6,5\}$ is an edge in $E_{3,4}$, then by Claim 1, there exists a unique distinguished edge containing $6$, written as $\{6,z''\}$, where $z''\in G_2$. Note that there are three distinguished edges of type $(0,2)$, and so far we have gotten two of them, i.e., $\{3,8\}$ and $\{9,11\}$. So the third one is just the one containing $6$. Since 5 is a c-element, the third one cannot contain $5$. Because of Claim 2, the third cannot be $\{6,8\}$ and $\{6,11\}$. Thus the third one would be $\{6,2\}$, but $\{6,2\}$ has been covered by the first block containing $4$, a contradiction.

W.l.o.g., set $x_3=3$, $x_4=6$ and $x_7=9$. There are three distinguished edges of type $(0,2)$, and so far we have gotten two of them, i.e., $\{6,8\}$ and $\{9,11\}$. Since 5 is a c-element, the third one cannot contain $5$. Since $d_3(6)=1$, the third one cannot contain $6$. Since $0$ is a vertex with degree one in $E_{3,4}$ and $\{0,1\}$ is an edge in $E_{3,4}$, then by Claim 1, there exists a unique distinguished edge containing $0$ which is of type $(0,1)$, and so the third cannot contain $0$. Because of Claim 2, the third one cannot be $\{3,8\}$, $\{3,11\}$, $\{9,2\}$ and $\{9,8\}$. Thus the third one would be $\{3,2\}$, but $\{3,2\}$ has been covered by the first block containing $4$, a contradiction.

(3.3) For $L=E_{3,5}$, since $4$ is a vertex with degree one in $E_{3,5}$ and $\{4,5\}$ is an edge in $E_{3,5}$, then by Claim 1, there exists a unique distinguished edge containing $4$, written as $\{4,z\}$, where $z\in \{2,8,11\}$. Since $8$ is also a vertex with degree one in $E_{3,5}$ and $\{8,6\}$ is an edge in $E_{3,5}$, then by Claim 1, there exists a unique distinguished edge containing $8$, written as $\{8,z'\}$, where $z'\in G_0$. Thus $z\neq 8$, and w.l.o.g., set $z=2$, i.e., $x_1=2$.

Since $0$ is a vertex with degree one in $E_{3,5}$ and $\{0,1\}$ is an edge in $E_{3,5}$, then by Claim 1, there exists a unique distinguished edge containing $0$, written as $\{0,z''\}$, where $z''\in G_1$. Thus $0\not\in\{x_4,x_5,x_7,x_8\}$, and so $x_2=0$. Then w.l.o.g., set $x_5=8$ and $x_8=11$. It follows that $\{x_3,x_4,x_7\}=\{3,6,9\}$. Similar arguments to those in (3.2), one can complete the proof. \qed

\section{$(K_4-e)$-IGDDs of type $g^{(n,h)}$}

A $G$-GDD of type $g^{n-h}(gh)^1$ is also called an {\em incomplete $G$-GDD} of type $g^{(n,h)}$, simply written as a $G$-IGDD of type $g^{(n,h)}$. The group of size $gh$ is called the {\em hole} of the IGDD. Filling in the hole of a $G$-IGDD of type $g^{(n,h)}$ by using a $G$-GDP of type $g^{h}$, one can obtain a $G$-GDP of type $g^{n}$.

\begin{Construction} \label{con:filling by MGDP}
Suppose that there exists a $G$-IGDD of type $g^{(n,h)}$. If there exists a $G$-GDP of type $g^{h}$ with leave having $L$ as its edge set, then there exists a $G$-GDP of type $g^{n}$ with leave having $L$ as its edge set.
\end{Construction}

In the rest of this section, we present recursive constructions and a partial existence result for $(K_4-e)$-IGDDs, which are used in a later section to produce certain $(K_4-e)$-MGDPs.

\subsection{Constructions from resolvable designs}

An $(H,G)$-design is called {\em resolvable} if its blocks can be partitioned into classes ${\cal P}_i$ (often referred to as {\em parallel classes}) such that every
vertex of $H$ appears in exactly one block of each ${\cal P}_i$. A resolvable $(H,G)$-design is sometimes also referred to as a {\em $G$-factorization} of $H$, and a class is called a {\em $G$-factor} of $H$. The case where $G=K_2$ (a single edge) is known as a 1-factorization. It is known that the complete $n$-partite equipartite graph $K_{n(g)}$, each partite set of size $g$, has a 1-factorization if and only if $gn$ is even (cf. \cite{y}).

\begin{Lemma} \label{lem:1-factorization}
There exist a $(K_{4}-e)$-IGDD of type $g^{((3h+2)/2,h)}$ for any positive integers $g$ and $h$ satisfying $g(h+2)\equiv 0\pmod{4}$ and $h\equiv 0\pmod{2}$.
\end{Lemma}

\proof Let $n=(3h+2)/2$. Let $G_i=\{nj+i:0\leq j\leq g-1\}$, $0\leq i\leq n-1$. We here construct a $(K_{4}-e)$-IGDD of type $g^{(n,h)}$ on $I_{gn}$ with the hole set $\bigcup_{i=0}^{h-1} G_i$ and groups $G_i$, $h\leq i\leq n-1$. Consider the complete $(n-h)$-partite equipartite graph $K_{(n-h)(g)}$ on vertex set $\bigcup_{i=h}^{n-1}G_i$ with partite sets $G_i$, $h\leq i\leq n-1$. Since $g(n-h)$ is even, it has a 1-factorization $\cal F$, which contains $(n-h-1)g=(hg)/2$ 1-factors. Let ${\cal F}=\{F_1,F_2,\ldots,F_{(hg)/2}\}$. Take any partition of $\bigcup_{i=0}^{h-1} G_i$ into $(gh)/2$ pairwise disjoint 2-subsets $\{x_l,y_l\}$, $1\leq l\leq (hg)/2$. Set
$${\cal B}=\{[f_1,f_2,x_l-y_l]:\{f_1,f_2\}\in F_l,1\leq l\leq (hg)/2\}.$$ It is readily checked that $\cal B$ is the set of blocks of the required $(K_{4}-e)$-IGDD of type $g^{(n,h)}$. \qed

Let $P_3$ be the path $[a;b,c]$ having vertex set $\{a,b,c\}$ and edges $\{a,b\}$, $\{a,c\}$. The following lemma was mentioned by Horton in \cite{h}, and was rediscovered by Milici and Tuza \cite{mt} recently.

\begin{Lemma} \label{lem:P_3}\emph{\cite{h,mt}}
Let $v\equiv 0\pmod{3}$ and $v\geq 9$. The union of any two edge-disjoint parallel classes of triangles of $K_v$ can be decomposed into three parallel classes of $P_3$.
\end{Lemma}

\begin{Lemma} \label{lem:K_4-e from P_3}
There exist a $(K_{4}-e)$-IGDD of type $3^{(n,3(n-1)/7)}$ for any $n\equiv 1\pmod{7}$ and $n\geq 8$.
\end{Lemma}

\proof Let $h=3(n-1)/7$. Let $G_i=\{3i+j:0\leq j\leq 2\}$, $0\leq i\leq n-1$. We here construct a $(K_{4}-e)$-IGDD of type $3^{(n,h)}$ on $I_{3n}$ with the hole set $I_{3h}=\bigcup_{i=0}^{h-1} G_i$ and groups $G_i$, $h\leq i\leq n-1$. It is known that there exists a resolvable $K_3$-GDD of type $1^u$, i.e., a Kirkman triple system of order $u$, for any $u\equiv 3\pmod{6}$, which is equivalent to a resolvable $K_3$-GDD of type $3^{u/3}$ (cf. \cite{StinsonBook}). Thus we can construct a resolvable $K_3$-GDD of type $3^{n-h}$ on $\bigcup_{i=h}^{n-1} G_i$ with groups $G_i$, $h\leq i\leq n-1$, which has $3(n-h-1)/2$ parallel classes. By Lemma \ref{lem:P_3}, it yields a resolvable $P_3$-GDD of type $3^{n-h}$, which has $9(n-h-1)/4=9(n-1)/7=3h$ parallel classes. Denote by ${\cal P}=\{P_0,P_1,\ldots,P_{3h-1}\}$ the set of its parallel classes. Set
$${\cal B}=\{[x,a_1,a_2-a_3]:[a_1;a_2,a_3]\in P_x,0\leq x\leq 3h-1\}.$$ It is readily checked that $\cal B$ is the set of blocks of the required $(K_{4}-e)$-IGDD of type $3^{(n,h)}$. \qed

\subsection{Inflation construction}

The {\em chromatic number} $\chi(G)$ of a graph $G$ is the least number of colors required
for a proper vertex coloring of $G$.

\begin{Construction} \label{con:inflation-general-1}
Suppose that there exists a $G$-IGDD of type $g^{(n,h)}$. Let $k=\chi(G)$. If there exists a $K_k$-GDD of type $t^k$, then there exists a $G$-IGDD of type $(gt)^{(n,h)}$.
\end{Construction}

\proof Let $V$ be a set of size $gn$ and can be partitioned into $n$ disjoint subsets $V_i$, $0\leq i\leq n-1$, each of size $g$. Assume that the given $G$-IGDD of type $g^{(n,h)}$ is constructed on $V$ with the hole set $\bigcup_{i=0}^{h-1} V_i$ and groups $V_i$, $h\leq i\leq n-1$. Denote by $\cal B$ the set of its blocks.

Since $\chi(G)=k$, for each $B\in{\cal B}$, we can give a proper vertex coloring of $B$ with $k$ colors. Such a vertex coloring is denoted by $\pi_B:B\rightarrow C=\{c_0,c_1,\ldots,c_{k-1}\}$.

Let $T$ be a set of size $t$. By assumption, we can construct a $K_k$-GDD of type $t^k$ on $C\times T$ with groups $\{c_j\}\times T$, $0\leq j\leq k-1$. Denote by $\cal A$ the set of its blocks. For each $A\in{\cal A}$, define a mapping $\rho_A:C\rightarrow T$ satisfying $\rho_A(c_j)=t_j$, where $A=\{(c_1,t_1),(c_2,t_2),\ldots,(c_k,t_k)\}$.

Now we construct a $G$-IGDD of type $(gt)^{(n,h)}$ on $V\times T$ with the hole set $\bigcup_{i=0}^{h-1} (V_i\times T)$ and groups $V_i\times T$, $h\leq i\leq n-1$. Let
$${\cal D}=\{\{(b,\rho_A(\pi_B(b))):b\in B\}:B\in{\cal B}, A\in{\cal A}\}.$$ It is readily checked that $\cal D$ forms the block set of the required $G$-IGDD of type $(gt)^{(n,h)}$. \qed

\begin{Construction} \label{con:inflation-1}
If there exists a $(K_{4}-e)$-IGDD of type $g^{(n,h)}$, then there exists a $(K_{4}-e)$-IGDD of type $(gt)^{(n,h)}$ for any $t\geq 1$.
\end{Construction}

\proof Obviously $\chi(K_{4}-e)=3$. It is well known that there exists a $K_3$-GDD of type $t^3$ for any $t\geq 1$, which is equivalent to a Latin square of order $t$ (cf. \cite{StinsonBook}). Then the conclusion follows from Construction \ref{con:inflation-general-1} straightforwardly.  \qed

The graph $K_4-e$ is equivalent to a $(K_{4}-e)$-IGDD of type $1^{(4,2)}$. So applying Construction \ref{con:inflation-1}, one can have

\begin{Lemma} \label{lem:n^2 2n^1}
There exists a $(K_{4}-e)$-IGDD of type $g^{(4,2)}$ for any $g\geq 1$.
\end{Lemma}

\begin{Lemma} \label{lem:g^(12,4)}
Let $n\in\{12,14\}$. There exists a $(K_{4}-e)$-IGDD of type $g^{(n,4)}$ for any $g\geq 1$.
\end{Lemma}

\proof Take a $(K_4-e)$-IGDD of type $1^{(n,4)}$ for $n\in\{12,14\}$, which exists by Examples 1.5 and 4.2 in \cite{hlss}. Then by Construction \ref{con:inflation-1}, we have a $(K_4-e)$-IGDD of type $g^{(n,4)}$ for any $g\geq 1$. \qed

\subsection{Filling constructions}

The following two constructions are variations of standard filling constructions for group divisible designs (cf. Lemmas 1.11 and 1.12 in \cite{chr}).

\begin{Construction} \label{con:filling by itselt}
Suppose that there exists a $G$-IGDD of type $g^{(n,h)}$. If there exists a $G$-IGDD of type $g^{(h,h_1)}$, then there exists a $G$-IGDD of type $g^{(n,h_1)}$.
\end{Construction}

\begin{Construction} \label{con:filling by GDD}
Suppose that there exists a $G$-GDD of type $(gm)^t$. If there exists a $G$-IGDD of type $g^{(m+h,h)}$, then there exist $G$-IGDDs of type $g^{(mt+h,m+h)}$ and $g^{(mt+h,h)}$.
\end{Construction}

\begin{Lemma} \label{lem:filling 1-type IGDD}
If there exists a $(K_{4}-e)$-IGDD of type $g^{(5+h,h)}$, then there exist $(K_{4}-e)$-IGDDs of type $g^{(n,5+h)}$ and $g^{(n,h)}$ for any $n\equiv h\pmod{5}$ and $n\geq 15+h$.
\end{Lemma}

\proof Let $t=(n-h)/5$. Then $t\geq 3$. By Theorem \ref{thm:K4-e-GDD}, there exists a $(K_4-e)$-GDD of type $(5g)^t$. Start from this GDD and apply Construction \ref{con:filling by GDD} with $m=5$ to complete the proof. \qed

\subsection{Direct constructions and the use of the recursions}

\begin{Lemma} \label{lem:IGDD g=2,3,4}
There exists a $(K_{4}-e)$-IGDD of type $g^{(5+h,h)}$ for $g\in\{2,3,4\}$ and $h\in\{2,3,4\}$.
\end{Lemma}
\proof Let $n=5+h$. For $(g,h)\in\{(2,3),(3,2),(3,4)\}$, we here construct a $(K_4-e)$-IGDD of type $g^{(n,h)}$ on $I_{gn}$ with the hole set $\{nj+i: 0\leq j\leq g-1, 0\leq i\leq h-1\}$ and groups $\{nj+i:0\leq j\leq (g-1)\}$, $h\leq i\leq n-1$. Let $\alpha_{g,h}=\prod_{i=0}^{n-1}(i\ i+n\ \ldots\ i+n(g-1))$ be a permutation on $I_{gn}$. Let $G_{g,h}$ be the group generated by $\alpha_{g,h}$. We only list base blocks. All other blocks are obtained by developing these base blocks under the action of $G_{g,h}$.
\begin{center}
\begin{tabular}{lllll}
$(g,h)=(2,3):$&
$[0,  4,  3-  5]$,&
$[0,  7,  6-  11]$,&
$[12,  14,  0-  1]$,&
$[13,  15,  0-  1]$,\\&
$[3,  5,  1-  10]$,&
$[4,  7,  1-  13]$,&
$[6,  11,  1-  5]$,&
$[2,  6,  3-  12]$,\\&
$[4,  15,  2-  11]$,&
$[2,  14,  5-  7]$.\\
$(g,h)=(3,2):$&
$[0,  3,  2-  4]$,&
$[0,  6,  5-  9]$,&
$[0,  12,  10-  11]$,&
$[0,  19,  13-  16]$,\\&
$[0,  20,  17-  18]$,&
$[1,  4,  2-  10]$,&
$[3,  12,  1-  16]$,&
$[1,  13,  5-  17]$,\\&
$[1,  18,  6-  9]$,&
$[11,  19,  1-  3]$,&
$[16,  20,  1-  10]$,&
$[2,  18,  12-  13]$.\\
$(g,h)=(3,4):$&
$[0,  5,  4-  6]$,&
$[0,  8,  7-  13]$,&
$[0,  16,  14-  15]$,&
$[0,  25,  17-  22]$,\\&
$[0,  26,  23-  24]$,&
$[1,  6,  4-  13]$,&
$[1,  15,  5-  7]$,&
$[1,  14,  8-  25]$,\\&
$[16,  23,  1-  11]$,&
$[17,  24,  1-  20]$,&
$[22,  26,  1-  11]$,&
$[4,  23,  2-  3]$,\\&
$[2,  22,  5-  16]$,&
$[8,  24,  2-  21]$,&
$[15,  25,  2-  3]$,&
$[3,  24,  5-  13]$,\\&
$[3,  22,  7-  8]$,&
$[3,  26,  14-  16]$.\\
\end{tabular}
\end{center}

For $(g,h)=(2,2)$, by Lemma \ref{lem:1-factorization}, there exists a $(K_4-e)$-IGDD of type $2^{(7,4)}$. Start from this IGDD and then apply Construction \ref{con:filling by itselt} with a $(K_4-e)$-IGDD of type $2^{(4,2)}$ (from Lemma \ref{lem:n^2 2n^1}), to obtain a $(K_4-e)$-IGDD of type $2^{(7,2)}$.

For $(g,h)=(2,4)$, we construct a $(K_4-e)$-IGDD of type $2^{(9,4)}$ on $I_{18}$ with the hole set $\{0,1,2,3,9,10,11,12\}$ and groups $\{i,i+9\}$, $4\leq i\leq 8$. Let $\alpha=(0\ 1\ 2\ 3\ 9\ 10\ 11\ 12)(4\ 5\ 6\ 7\ 13\ 14\ 15$ $16)$ be a permutation on $I_{18}$. Let $G$ be the group generated by $\alpha$. Base blocks are listed below, and all other blocks are obtained by developing these base blocks under the action of $G$.
\begin{center}
\begin{tabular}{lll}
$[0,  5,  4-  7]$,&
$[6,  14,  0-  11]$,&
$[0,  15,  8-  17]$.\\
\end{tabular}
\end{center}

For $(g,h)=(3,3)$, the conclusion follows from Lemma \ref{lem:K_4-e from P_3}. For $g=4$ and $h\in\{2,3,4\}$, start from a $(K_{4}-e)$-IGDD of type $2^{(n,h)}$, and then apply Construction \ref{con:inflation-1} with $t=2$ to obtain a $(K_{4}-e)$-IGDD of type $4^{(n,h)}$. \qed

\begin{Lemma} \label{lem:IGDD-2^(14,7)}
There exists a $(K_{4}-e)$-IGDD of type $2^{(14,7)}$.
\end{Lemma}

\proof Let $G_i=\{i,i+14\}$, $0\leq i\leq 13$. We construct the required design on $I_{28}$ with the hole set $\bigcup_{i=0}^6 G_i$ and groups $G_i$, $7\leq i\leq 13$. Let $\alpha=(0\ 1\ 2\ 3\ 4\ 5\ 6\ 14\ 15\ 16\ 17\ 18\ 19\ 20)$ $(7\ 8\ 9\ 10\ 11\ 12\ 13\ 21\ 22\ 23\ 24\ 25\ 26\ 27)$ be a permutation on $I_{28}$. Let $G$ be the group generated by $\alpha$. We only list base blocks below.

\begin{center}
\begin{tabular}{llllll}
$[0,  8,  7-  10]$,&
$[0,  13,  9-  23]$,&
$[11,  24,  0-  20]$,&
$[21,  26,  0-  20]$.\\
\end{tabular}
\end{center}
All other blocks are obtained by developing these base blocks under the action of $G$. \qed

\begin{Lemma} \label{lem:11^(8,3)}
There exists a $(K_{4}-e)$-IGDD of type $11^{(8,3)}$.
\end{Lemma}
\proof Let $G_i=\{8j+i: 0\leq j\leq 10\}$, $0\leq i\leq 7$. We construct the required design on $I_{88}$ with the hole set $\bigcup_{i=0}^2 G_i$ and groups $G_i$, $3\leq i\leq 7$. Let $\alpha=(0\ 8\ 16 \cdots 80)(1\ 9\ 17 \cdots 81)(2\ 10\ 18 \cdots 82)$ $(3\ 4 \cdots 7\ 11\ 12 \cdots 15\cdots 83\ 84 \cdots 87)$ be a permutation on $I_{88}$. Let $G$ be the group generated by $\alpha$. Only base blocks are listed below.

\begin{center}
\begin{tabular}{llllll}
$[0,  4,  3-  6]$,&
$[0,  11,  5-  15]$,&
$[0,  19,  7-  29]$,&
$[13,  31,  0-  1]$,&
$[1,  15,  3-  4]$,&
$[1,  27,  5-  53]$,\\
$[1,  45,  6-  12]$,&
$[2,  28,  3-  63]$,&
$[4,  35,  2-  55]$,&
$[5,  51,  2-  74]$,&
$[6,  47,  2-  67]$.
\end{tabular}
\end{center}
All other blocks are obtained by developing these base blocks under the action of $G$. \qed

\section{$(K_{4}-e)$-IGDDs of type $(g,h)^n$}

We now introduce another type of IGDDs, which are frequently used later. Let $S$ be a subset of vertex set of the complete $n$-partite equipartite graph $K_{n(g)}$ with partite sets $G_i$ (called {\em groups}), $1\leq i\leq n$, satisfying $|G_i\cap S|=h$. Let $K[S]$ be the subgraph of $K_{n(g)}$ induced by $S$. A $(K_{n(g)}\setminus K[S],G)$-design is also refer to as an {\em incomplete group divisible design} and $S$ is called the {\em hole}. Such a design is denoted by a $G$-IGDD of type $(g,h)^n$.

\begin{Example} \label{eg:4^3-IGDD}
We here explicitly construct a $(K_{4}-e)$-IGDD of type $(4,1)^3$ on $I_{12}$ with groups $\{i,i+3,i+6,i+9\}$, $0\leq i\leq 2$, and the hole set $S=\{0,1,2\}$. All its blocks are:
\begin{center}
\begin{tabular}{llllll}
$[5,  4,  0-  3]$,&
$[7,  8,  0-  6]$,&
$[10,  11,  0-  9]$,&
$[5,  6,  1-  10]$,&
$[3,  11,  1-  7]$,&
$[8,  9,  1-  4]$,\\
$[7,  9,  2-  5]$,&
$[3,  10,  2-  8]$,&
$[4,  6,  2-  11]$.\\
\end{tabular}
\end{center}
\end{Example}

Filling in the hole of a $G$-IGDD of type $(g,h)^n$ by using a $G$-GDP of type $h^{n}$, one can obtain a $G$-GDP of type $g^{n}$.

\begin{Construction} \label{con:fillling-2}
Suppose that there exists a $G$-IGDD of type $(g,h)^n$. If there exists a $G$-GDP of type $h^{n}$ with leave having $L$ as its edge set, then there exists a  $G$-GDP of type $g^{n}$ with leave having $L$ as its edge set.
\end{Construction}

\begin{Example}\label{eg:4^3}
Take a $(K_{4}-e)$-IGDD of type $(4,1)^3$ from Example $\ref{eg:4^3-IGDD}$. Then applying Construction $\ref{con:fillling-2}$ with a trivial $(K_{4}-e)$-MGDP of type $1^3$, one can obtain a $(K_{4}-e)$-MGDP of type $4^3$ with leave $E_{3,1}$.
\end{Example}

\begin{Example} \label{eg:IGDD-(4,1)^4}
There exist a $(K_{4}-e)$-IGDD of type $(4,1)^4$ and a $(K_{4}-e)$-MGDP of type $4^4$ with leave $E_{1}$.
\end{Example}

\proof We construct a $(K_{4}-e)$-IGDD of type $(4,1)^4$ below on $I_{16}$ with groups $\{i,i+4,i+8,i+12\}$, $0\leq i\leq 3$, and the hole set $\{0,1,14,15\}$:

\begin{center}
\begin{tabular}{llllll}
$[0,  3,  2-  5]$,&
$[0,  7,  6-  9]$,&
$[0,  11,  10-  13]$,&
$[5,  15,  10-  12]$,&
$[1,  4,  2-  3]$,&
$[1,  8,  6-  7]$,\\
$[2,  7,  5-  12]$,&
$[2,  9,  8-  11]$,&
$[1,  12,  10-  11]$,&
$[13,  15,  2-  8]$,&
$[6,  12,  3-  13]$,&
$[3,  10,  8-  9]$,\\
$[13,  14,  3-  7]$,&
$[5,  11,  4-  6]$,&
$[8,  14,  5-  11]$,&
$[4,  10,  7-  13]$,&
$[9,  14,  4-  12]$,&
$[6,  15,  4-  9]$.\\
\end{tabular}
\end{center}
Then apply Construction \ref{con:fillling-2} with a $(K_{4}-e)$-MGDP of type $1^4$, which has only one block, to obtain a $(K_{4}-e)$-MGDP of type $4^4$ with leave $E_{1}$. \qed

In the rest of this section, we present recursive constructions and a partial existence result for such kind of $(K_4-e)$-IGDDs.

\subsection{Constructions from $2$-frames}

Recall that a $K_2$-GDD of type $g^n$ consists of all edges of a complete $n$-partite equipartite graph. A $K_2$-GDD of type $g^n$ defined on $I_{gn}$ with groups $G_1,G_2,\ldots,G_n$ is called a {\em $2$-frame} if its blocks can be partitioned into $gn$ classes ${\cal P}_{i,j}$, $1\leq i\leq n$, $1\leq j\leq g$, (often referred to as {\em partial parallel classes}) such that every vertex of $I_{gn}\setminus G_i$, $1\leq i\leq n$, appears in exactly one block of ${\cal P}_{i,j}$ for any $1\leq j\leq g$.

\begin{Lemma} \label{lem:2-frame}\emph{\cite{rs}}
There exists a $2$-frame of type $g^n$ if and only if $g(n-1)\equiv 0 \pmod{2}$ and $n\geq 3$.
\end{Lemma}

\begin{Lemma} \label{lem:IGDD-(3,2)^n}
There exists a $(K_{4}-e)$-IGDD of type $(3g,2g)^{n}$ for any $g(n-1)\equiv 0 \pmod{2}$ and $n\geq 3$.
\end{Lemma}

\proof We here construct a $(K_{4}-e)$-IGDD of type $(3g,2g)^{n}$ on $I_{3gn}$ with groups $G_i=\{nj+i:0\leq j\leq 3g-1\}$, $0\leq i\leq n-1$, and the hole set $I_{3gn}\setminus I_{gn}$. First we construct a $2$-frame of type $g^n$ on $I_{gn}$ with groups $\{nj+i:0\leq j\leq g-1\}$, $0\leq i\leq n-1$, which exists for $g(n-1)\equiv 0 \pmod{2}$ and $n\geq 3$ by Lemma \ref{lem:2-frame}. It has $gn$ partial parallel classes ${\cal P}_{i,j}$, $0\leq i\leq n-1$, $0\leq j\leq g-1$. For each $0\leq i\leq n-1$, take any partition of $\{nj+i:g\leq j\leq 3g-1\}$ into $g$ pairwise disjoint 2-subsets $\{x^{(i)}_l,y^{(i)}_l\}$, $0\leq l\leq g-1$. Set
$${\cal B}=\{[p_1,p_2,x^{(i)}_j-y^{(i)}_j]:\{p_1,p_2\}\in {\cal P}_{i,j},0\leq i\leq n-1,0\leq j\leq g-1\}.$$ It is readily checked that $\cal B$ is the set of blocks of a $(K_{4}-e)$-IGDD of type $(3g,2g)^{n}$. \qed

\subsection{Constructions from pairwise balanced designs}

Let $L=\{k_1,k_2,\ldots,k_r\}$ be a set of positive integers and ${\cal G}=\{K_{k_1},K_{k_2},\ldots,K_{k_r}\}$ be a set of cliques. A $(K_n, {\cal G})$-design is said to be a {\em pairwise balanced design}, denoted by an $(n,L,1)$-PBD.

\begin{Construction}\label{con:PBD}\emph{\cite{hcf}}
Suppose there exists an $(n,L,1)$-PBD. If for each $k\in L$, there exists a
$G$-IGDD of type $(g,h)^k$, then there exists a $G$-IGDD of type $(g,h)^n$.
\end{Construction}

\begin{Lemma}\label{lem:PBD(3,4,5)}\emph{\cite{abg}}
There exists a $(n,\{3,4,5\},1)$-PBD for any $n\geq 3$ and $n\neq 6,8$.
\end{Lemma}

\subsection{Inflation construction}

By similar arguments to those in Construction \ref{con:inflation-general-1}, one can obtain

\begin{Construction} \label{con:inflation-general-2}
Suppose that there exists a $G$-IGDD of type $(g,h)^{n}$. Let $k=\chi(G)$. If there exists a $K_k$-GDD of type $t^k$, then there exists a $G$-IGDD of type $(gt,ht)^{n}$.
\end{Construction}

Applying Construction \ref{con:inflation-general-2} with a $K_3$-GDD of type $t^3$, which is equivalent to a Latin square of order $t$ (cf. \cite{StinsonBook}) and exists for any $t\geq 1$, one can have

\begin{Construction} \label{con:inflation-3}
If there exists a $(K_{4}-e)$-IGDD of type $(g,h)^{n}$, then there exists a $(K_{4}-e)$-IGDD of type $(gt,ht)^{n}$ for any $t\geq 1$.
\end{Construction}

\subsection{Filling constructions}

The following construction is a variation of Construction 2.5 in \cite{w}.

\begin{Construction}\label{con:filling IGDD-1}
Suppose that there exists a $G$-IGDD of type $(g,h)^{n}$. If there exists a $G$-IGDD of type $(h,r)^{n}$, then there exists a $G$-IGDD of type $(g,r)^{n}$.
\end{Construction}

Let $S_1,S_2,\ldots,S_{u+1}$ be disjoint subsets of vertex set of the complete $n$-partite equipartite graph $K_{n(hu+w)}$ with partite sets $G_i$ (called {\em groups}), $1\leq i\leq n$, satisfying $|S_j|=hn$ for $1\leq j\leq u$, $|S_{u+1}|=wn$, and $|G_i\cap S_j|=h$ for any $1\leq i\leq n$ and $1\leq j\leq u$, $|G_i\cap S_{u+1}|=w$ for any $1\leq i\leq n$. Let $K[S_j]$ be the subgraph of $K_{n(hu+w)}$ induced by $S_j$. A $(K_{n(hu+w)}\setminus (\bigcup_{j=1}^{u+1} K[S_j]),G)$-design is referred to as a {\em holey group divisible design} and each $S_j$ is called a {\em hole}. Such a design is denoted by a $G$-HGDD of type $(n,h^u w^1)$.

The following construction is a variation of Construction 2.9 in \cite{w}.

\begin{Construction}\label{con:filling HGDD}
Suppose that there exists a $G$-HGDD of type $(n,h^{u} w^1)$.
\begin{enumerate}
\item[$(1)$] If there exists a $G$-IGDD of type $(h+r,r)^{n}$, then there exists a $G$-IGDD of type $(hu+w+r,w+r)^{n}$.
\item[$(2)$] If there exist $G$-IGDDs of type $(h+r,r)^{n}$ and $(w+r,r)^{n}$, then there exist $G$-IGDDs of type $(hu+w+r,h+r)^{n}$ and $(hu+w+r,r)^{n}$.
\end{enumerate}
\end{Construction}

\begin{Construction}\label{con:from HGDD}
Suppose that there exists a $G$-GDD of type $h^u w^1$. Let $k=\chi(G)$. If there exists a $K_k$-HGDD of type $(k,1^{n})$, then there exists a $G$-HGDD of type $(n,h^u w^1)$.
\end{Construction}

\proof Let $V$ be a set of size $hu+w$ and can be partitioned into $u+1$ disjoint subsets $V_0,V_1,\ldots,V_{u-1},W$, where $|V_i|=h$ for $0\leq i\leq u-1$ and $|W|=w$. Assume that the given $G$-GDD of type $h^u w^1$ is constructed on $V$ with groups $V_0,V_1,\ldots,V_{u-1},W$. Denote by $\cal B$ the set of its blocks.

Since $\chi(G)=k$, for each $B\in{\cal B}$, we can give a proper vertex coloring of $B$ with $k$ colors. Such a vertex coloring is denoted by $\pi_B:B\rightarrow C=\{c_0,c_1,\ldots,c_{k-1}\}$.

Let $N$ be a set of size $n$. By assumption, we can construct a $K_k$-HGDD of type $(k,1^{n})$ on $C\times N$ with groups $\{c_j\}\times N$, $0\leq j\leq k-1$, and holes $\{x\}\times C$, $x\in N$. Denote by $\cal A$ the set of its blocks. For each $A\in{\cal A}$, define a mapping $\rho_A:C\rightarrow N$ satisfying $\rho_A(c_j)=n_j$, where $A=\{(c_1,n_1),(c_2,n_2),\ldots,(c_k,n_k)\}$.

Now we construct a $G$-HGDD of type $(n,h^u w^1)$ on $N\times V$ with groups $\{x\}\times V$, $x\in N$, and holes $N\times W$, $N\times V_i$, $0\leq i\leq u-1$. Let
$${\cal D}=\{\{(\rho_A(\pi_B(b)),b):b\in B\}:B\in{\cal B}, A\in{\cal A}\}.$$ It is readily checked that $\cal D$ forms the block set of the required $G$-HGDD of type $(n,h^u w^1)$. \qed

\begin{Construction} \label{con:from HGDD-d}
If there exists a $(K_{4}-e)$-GDD of type $h^u w^1$, then there exists a $(K_{4}-e)$-HGDD of type $(n,h^u w^1)$ for any $n\geq 3$.
\end{Construction}

\proof Begin with the given $(K_{4}-e)$-GDD of type $h^u w^1$, and then apply Construction \ref{con:from HGDD} with a $K_3$-HGDD of type $(3,1^{n})$, which is equivalent to an idempotent Latin square of order $n$ and exists for any $n\geq 3$ (cf. \cite{StinsonBook}), to obtain a $(K_{4}-e)$-HGDD of type $(n,h^u w^1)$. \qed

\begin{Lemma}\label{lem:HGDD-(n,5^u)}
\begin{enumerate}
\item[$(1)$] There exists a $(K_{4}-e)$-HGDD of type $(n,h^u)$ for any $n\geq 3$, $u\geq 3$, $hu(u-1)\equiv 0\pmod{5}$ and $(h,u)\neq(1,5)$.
\item[$(2)$] There exists a $(K_{4}-e)$-HGDD of type $(n,h^2 (2h)^1)$ for any $n\geq 3$ and $h\geq 1$.
\end{enumerate}
\end{Lemma}

\proof By Theorem \ref{thm:K4-e-GDD}, there exists a $(K_4-e)$-GDD of type $h^u$ for any $u\geq 3$, $hu(u-1)\equiv 0\pmod{5}$ and $(h,u)\neq(1,5)$. By Lemma \ref{lem:n^2 2n^1}, there exists a $(K_{4}-e)$-IGDD of type $h^{(4,2)}$, i.e., a $(K_{4}-e)$-GDD of type $h^2 (2h)^1$, for any $h\geq 1$. Start from these GDDs and apply Construction \ref{con:from HGDD-d} to complete the proof. \qed

\begin{Lemma} \label{lem:filling 2-type IGDD}
If there exists a $(K_{4}-e)$-IGDD of type $(5+r,r)^{n}$, then there exist $(K_{4}-e)$-IGDDs of type $(g,5+r)^{n}$ and $(g,r)^{n}$ for any $g\equiv r\pmod{5}$ and $g\geq 15+r$.
\end{Lemma}

\proof Let $u=(g-r)/5$. Then $u\geq 3$. Begin with a $(K_{4}-e)$-HGDD of type $(n,5^u)$ (from Lemma \ref{lem:HGDD-(n,5^u)}(1)). Then apply Construction \ref{con:filling HGDD} with the given $(K_{4}-e)$-IGDD of type $(5+r,r)^{n}$ to obtain $(K_{4}-e)$-IGDDs of type $(5u+r,5+r)^{n}$ and $(5u+r,r)^{n}$. \qed

\subsection{Direct constructions and the use of the recursions}

\begin{Lemma}\label{lem:IGDD-(4,1)^5)}
There exists a $(K_{4}-e)$-IGDD of type $(4,1)^n$ for any $n\geq 3$ and $n\neq 6$.
\end{Lemma}

\proof For $n=3,4$, the conclusion follows from Examples \ref{eg:4^3-IGDD} and \ref{eg:IGDD-(4,1)^4}. For $n=5,8$, we construct the required designs on $I_{4n}$ with groups $\{nj+i:0\leq j\leq 3\}$, $0\leq i\leq n-1$, and the hole set $I_n$. Let $\alpha_5=(0\ 1\ 2\ 3\ 4)(5\ 6\ 7\ 8\ 9)(10\ 11\ 12\ 13\ 14)(15\ 16\ 17\ 18\ 19)$ be a permutation on $I_{20}$. Let $\alpha_8=(0\ 2\ 4\ 6)(1\ 3\ 5\ 7)(8\ 10\ 12\ 14)(9\ 11\ 13\ 15)(16\ 18\ 20\ 22)(17\ 19\ 21\ 23)$ $(24\ 26\ 28\ 30)(25\ 27\ 29\ 31)$ be a permutation on $I_{32}$. Let $G_n$ be the group generated by $\alpha_n$, $n=5,8$. The base blocks are listed below.
\begin{center}{ \tabcolsep 0.045in
\begin{tabular}{llllll}
$n=5$:& $[ 0, 6, 7-8 ]$,& $[ 0, 9, 11-12 ]$,& $[ 15, 16, 2-4 ]$,& $[ 10, 11, 2-18 ]$,& $[ 15, 17, 8-9 ]$,\\& $[ 10, 12, 6-16 ]$.\\
$n=8$:& $[8, 12, 17- 21]^*$,&
$[9, 13, 16- 20]^*$,&
$[16, 20, 25- 29]^*$,&
$[17, 21, 24- 28]^*$,&
$[24, 28, 9- 13]^*$,\\&
$[25, 29, 8- 12]^*$,&
$[0,  10,  9-  11]$,&
$[0,  14,  12-  17]$,&
$[0,  15,  13-  19]$,&
$[0,  20,  18-  21]$,\\&
$[0,  25,  22-  23]$,&
$[0,  27,  26-  30]$,&
$[0,  31,  28-  29]$,&
$[1,  11,  8-  16]$,&
$[1,  20,  10-  14]$,\\&
$[12,  26,  1-  16]$,&
$[13,  31,  1-  8]$,&
$[1,  24,  15-  19]$,&
$[1,  21,  18-  30]$,&
$[1,  29,  22-  23]$,\\&
$[27,  28,  1-  15]$,&
$[23,  27,  8-  13]$,&
$[26,  28,  8-  22]$,&
$[18,  23,  9-  17]$.\\
\end{tabular}}
\end{center}
All other blocks are obtained by developing these base blocks under the action of $G_n$. Note that each base block marked $^*$ for $n=8$ only generates $2$ different blocks under the action of $G_8$.

For $n\geq 7$ and $n\neq 8$, start from an $(n,\{3,4,5\},1)$-PBD (from Lemma \ref{lem:PBD(3,4,5)}), and then apply Construction \ref{con:PBD} with a $(K_4-e)$-IGDD of type $(4,1)^k$, $k\in\{3,4,5\}$. \qed

\begin{Lemma}\label{lem:(7,3)^n}
There exists a $(K_4-e)$-IGDD of type $(7,3)^n$ for any $n\geq 3$ and $n\neq 6$.
\end{Lemma}

\proof For $n\in\{3,4,5,8\}$, we construct a $(K_4-e)$-IGDD of type $(7,3)^n$ on $I_{7n}$ with groups $\{nj+i: 0\leq j\leq 6\}$, $0\leq i\leq n-1$, and the hole set $I_{3n}$. Let $\alpha_{n}=(0\ 1 \cdots n-1)(n\ n+1 \cdots 2n-1)(2n\ 2n+1 \cdots 3n-1)(3n\ 3n+1 \cdots 7n-1)$ be a permutation on $I_{7n}$. Let $G_{n}$ be the group generated by $\alpha_{n}$. Note that $|G_n|=4n$. Only base blocks are listed below.

\begin{center}
\begin{tabular}{lllll}
$n=3$:&
$[10,  11,  0-  15]$,&
$[10,  20,  3-  6]$.\\
$n=4$:&
$[0,  14,  13-  19]$,&
$[4,  15,  13-  18]$,&
$[8,  19,  13-  26]$.\\
$n=5$:&
$[16,17,0-3]$,&
$[16,18,5-9]$,&
$[16,22,10-13]$,&
$[15,26,18-19]$.\\
$n=8$:&
$[0,26,25-28]$,&
$[27,30,0-8]$,&
$[29,39,0-8]$,&
$[8,34,25-28]$,\\&
$[16,29,25-42]$,&
$[27,38,16-45]$,&
$[28,55,16-40]$.\\
\end{tabular}
\end{center}
All other blocks are obtained by developing these base blocks under the action of $G_{n}$.

For $n\geq 7$ and $n\neq 8$, start from an $(n,\{3,4,5\},1)$-PBD (from Lemma \ref{lem:PBD(3,4,5)}), and then apply Construction \ref{con:PBD} with a $(K_4-e)$-IGDD of type $(7,3)^k$, $k\in\{3,4,5\}$, to obtain a $(K_4-e)$-IGDD of type $(7,3)^n$. \qed

\begin{Lemma}\label{lem:(5+r,r)^8}
There exists a $(K_4-e)$-IGDD of type $(5+r,r)^8$ for $r\in\{2,3,4\}$.
\end{Lemma}

\proof We construct a $(K_4-e)$-IGDD of type $(5+r,r)^8$ on
$I_{40+8r}$ with groups $\{8j+i: 0\leq j\leq 4+r\}$, $0\leq
i\leq 7$, and the hole set $I_{8r}$. Let $\alpha_{r}=\prod_{j=0}^{4+r}(8j\ 8j+1\ 8j+2\ 8j+3\ 8j+4\ 8j+5\ 8j+6\ 8j+7)$ be a permutation on $I_{40+8r}$. Let $G_{r}$ be the group generated by $\alpha_{r}$. Only base blocks are listed below.

\begin{center}{\tabcolsep 0.02in
\begin{tabular}{llllll}
$r=2$:& $[16,20,26-30]^*$, & $[26,30,32-36]^*$, & $[32,36,42-46]^*$, & $[42,46,48-52]^*$, & $[48,52,18-22]^*$, \\& $[37,52,0-8]$,& $[0,41,38-46]$,&
$[8,54,17-35]$,&
$[8,49,18-51]$,&
$[8,21,19-30]$,\\&
$[8,31,20-26]$,&
$[8,36,22-34]$,&
$[8,27,23-25]$,&
$[8,33,28-29]$,&
$[8,55,38-46]$,\\&
$[39,43,8-22]$,&
$[41,42,8-13]$,&
$[8,50,47-53]$,&
$[16,46,19-39]$,&
$[31,44,16-49]$,\\&
$[34,47,16-25]$,&
$[16,49,35-37]$,&
$[0,18,17-31]$,&
$[0,39,19-28]$,&
$[0,45,20-25]$,\\&
$[0,47,21-29]$,&
$[0,50,22-23]$,&
$[26,43,0-45]$,&
$[27,34,0-52]$,&
$[30,53,0-31]$,\\&
$[0,36,33-35]$,&
$[42,49,0-6]$,&
$[54,55,0-26]$,&
$[47,51,24-38]$.\\

$r=3$:& $[24,28,34-38]^*$, & $[34,38,40-44]^*$, & $[40,44,50-54]^*$, & $[50,54,56-60]^*$, & $[56,60,26-30]^*$, \\&
$[59,62,0-42]$,&
$[8,50,25-27]$,&
$[26,38,8-52]$,&
$[28,42,8-13]$,&
$[29,44,8-62]$,\\&
$[8,51,30-39]$,&
$[33,34,8-13]$,&
$[8,52,35-54]$,&
$[41,60,8-10]$,&
$[43,55,8-17]$,\\&
$[46,49,8-34]$,&
$[53,62,8-16]$,&
$[57,59,8-12]$,&
$[16,27,25-26]$,&
$[16,31,28-50]$,\\&
$[16,41,29-30]$,&
$[33,35,16-23]$,&
$[24,60,35-59]$,&
$[37,55,16-34]$,&
$[16,57,38-39]$,\\&
$[16,58,43-44]$,&
$[45,52,16-24]$,&
$[46,59,16-23]$,&
$[0,34,25-27]$,&
$[0,39,26-42]$,\\&
$[0,46,28-33]$,&
$[0,58,29-35]$,&
$[0,61,30-54]$,&
$[31,53,0-40]$,&
$[0,57,36-37]$,\\&
$[38,45,0-47]$,&
$[0,50,41-63]$,&
$[0,44,43-47]$,&
$[49,60,0-38]$,&
$[0,52,51-55]$,\\&
$[49,61,16-22]$.\\

$r=4$:& $[32,36,42-46]^*$, & $[42,46,48-52]^*$, & $[48,52,58-62]^*$, & $[58,62,64-68]^*$, & $[64,68,34-38]^*$, \\&
$[24,60,35-37]$,&
$[24,71,36-49]$,&
$[41,59,24-28]$,&
$[24,70,42-57]$,&
$[0,42,33-37]$,\\&
$[0,71,34-35]$,&
$[0,57,36-66]$,&
$[0,52,38-61]$,&
$[0,46,39-63]$,&
$[0,70,41-69]$,\\&
$[0,55,43-44]$,&
$[0,50,45-51]$,&
$[0,58,47-54]$,&
$[49,68,0-21]$,&
$[0,65,53-67]$,\\&
$[0,62,59-60]$,&
$[33,34,8-11]$,&
$[8,66,35-47]$,&
$[8,62,36-69]$,&
$[8,57,37-50]$,\\&
$[41,43,8-29]$,&
$[8,65,42-68]$,&
$[16,67,42-45]$,&
$[43,69,16-24]$,&
$[16,59,46-70]$,\\&
$[47,62,16-24]$,&
$[51,68,16-24]$,&
$[16,50,53-65]$,&
$[44,53,8-15]$,&
$[8,58,46-55]$,\\&
$[49,51,8-31]$,&
$[52,70,8-29]$,&
$[59,71,8-12]$,&
$[60,61,8-16]$,&
$[33,66,16-24]$,\\&
$[16,36,34-63]$,&
$[35,57,16-28]$,&
$[37,54,16-23]$,&
$[39,49,16-36]$,&
$[16,44,41-58]$,\\&
$[24,54,34-67]$,&
$[24,53,38-58]$,&
$[32,44,43-51]$.
\end{tabular}}
\end{center}
All other blocks are obtained by developing these base blocks under the action of $G_{r}$. Note that each base block marked $^*$ only generates $4$ different blocks under the action of $G_r$. \qed

\begin{Lemma}
\label{lem:IGDD-(6,1)^n,(7,2)^n,(8,3)^n,(9,4)^n}
Let $r\in\{1,2,3,4\}$. There exists a $(K_4-e)$-IGDD of type $(5+r,r)^n$ for any $n\geq 3$ and $n\neq 6$.
\end{Lemma}

\proof For $r=1$, $n\geq 3$ and $n\neq 6$, take a $(K_4-e)$-IGDD of type $(6,4)^n$ (from Lemma \ref{lem:IGDD-(3,2)^n}), and then apply Construction \ref{con:filling IGDD-1} with a $(K_4-e)$-IGDD of type $(4,1)^n$ (from Lemma \ref{lem:IGDD-(4,1)^5)}) to obtain a $(K_4-e)$-IGDD of type $(6,1)^n$.

For $r=2$ and $n\in\{3,5\}$, take a $(K_4-e)$-IGDD of type $(7,3)^n$ (from Lemma \ref{lem:(7,3)^n}), and then apply Construction \ref{con:filling IGDD-1} with a $(K_4-e)$-IGDD of type $(3,2)^n$ (from Lemma \ref{lem:IGDD-(3,2)^n}) to obtain a $(K_4-e)$-IGDD of type $(7,2)^n$.

For $r=3$ and $n\in\{3,5\}$, start from a $(K_4-e)$-HGDD of type $(n,1^6)$ (from Lemma \ref{lem:HGDD-(n,5^u)}(1)), and then apply Construction \ref{con:filling HGDD} with a $(K_4-e)$-IGDD of type $(3,2)^n$ (from Lemma \ref{lem:IGDD-(3,2)^n}) to obtain a $(K_4-e)$-IGDD of type $(8,3)^n$.

For $r=4$ and $n\in\{3,5\}$, take a $(K_4-e)$-IGDD of type $(9,6)^n$ (from Lemma \ref{lem:IGDD-(3,2)^n}), and then apply Construction \ref{con:filling IGDD-1} with a $(K_4-e)$-IGDD of type $(6,4)^n$ (from Lemma \ref{lem:IGDD-(3,2)^n}) to obtain a $(K_4-e)$-IGDD of type $(9,4)^n$.

For $r\in\{2,3,4\}$ and $n=4$, we construct a $(K_4-e)$-IGDD of type $(5+r,r)^4$ on
$I_{20+4r}$ with groups $\{4j+i: 0\leq j\leq 4+r\}$, $0\leq
i\leq 3$, and the hole set $\{20,21,\ldots,19+4r\}$. Let $\alpha_{r}=\prod_{j=0}^{4+r}(4j\ 4j+1\ 4j+2)$ be a permutation on $I_{20+4r}$. Let $G_{r}$ be the group generated by $\alpha_{r}$. Only base blocks are listed below.

\begin{center}{\tabcolsep 0.04in
\begin{tabular}{llllll}
$r=2$:&
$[0,14,1-3]$,&
$[0,7,5-18]$,&
$[0,17,6-27]$,&
$[3,17,4-24]$,&
$[19,25,4-0]$,\\&
$[9,23,0-4]$,&
$[10,21,0-3]$,&
$[11,22,0-17]$,&
$[15,26,0-12]$,&
$[4,22,5-15]$,\\&
$[4,27,10-14]$,&
$[4,13,11-26]$,&
$[8,25,7-11]$,&
$[12,22,7-13]$,&
$[8,14,9-19]$,\\&
$[8,17,15-26]$,&
$[18,21,8-19]$,&
$[12,18,17-23]$.\\

$r=3$:&
$[18,25,4-7]$,&
$[11,26,8-17]$,&
$[13,27,8-18]$,&
$[8,25,14-19]$,&
$[0,22,11-13]$,\\&
$[3,21,8-14]$,&
$[18,23,0-9]$,&
$[0,26,1-3]$,&
$[0,6,5-27]$,&
$[0,14,7-31]$,\\&
$[0,15,9-29]$,&
$[0,17,10-30]$,&
$[19,21,0-6]$,&
$[3,29,4-18]$,&
$[4,21,7-10]$,\\&
$[4,31,9-17]$,&
$[11,30,4-12]$,&
$[4,14,13-23]$,&
$[15,26,4-13]$,&
$[8,29,7-10]$,\\&
$[12,19,18-29]$,&
$[16,21,15-18]$.\\

$r=4$:&
$[23, 1, 16-10]$,&
$[4, 14, 23-19]$,&
$[17, 19, 24-32]$,&
$[1, 14, 24-11]$,&
$[24, 5, 2-15]$,\\&
$[28, 3, 2-14]$,&
$[5, 11, 28-32]$,&
$[31, 9, 16-4]$,&
$[20, 18, 17-11]$,&
$[20, 1, 19-2]$,\\&
$[3, 6, 20-32]$,&
$[20, 14, 5-15]$,&
$[7, 9, 20-34]$,&
$[10, 13, 20-28]$,&
$[18, 3, 24-9]$,\\&
$[24, 7, 6-13]$,&
$[24, 10, 9-11]$,&
$[27, 2, 16-4]$,&
$[8, 13, 27-34]$,&
$[17, 6, 28-4]$,\\&
$[28, 7, 18-1]$,&
$[28, 9, 19-15]$,&
$[0, 14, 31-33]$,&
$[32, 15, 18-1]$,&
$[16, 13, 35-14]$,\\&
$[35, 10, 0-4]$.\\
\end{tabular}}
\end{center}
All other blocks are obtained by developing these base blocks under the action of $G_{r}$.

For $r\in\{2,3,4\}$ and $n=8$, the conclusion follows from Lemma \ref{lem:(5+r,r)^8}. For $r\in\{2,3,4\}$, $n\geq 7$ and $n\neq 8$, start from an $(n,\{3,4,5\},1)$-PBD (from Lemma \ref{lem:PBD(3,4,5)}), and then apply Construction \ref{con:PBD} with a $(K_4-e)$-IGDD of type $(5+r,r)^k$, $k\in\{3,4,5\}$, to obtain a $(K_4-e)$-IGDD of type $(5+r,r)^n$. \qed

\begin{Lemma}\label{lem:(11,1)^n}
There exists a $(K_4-e)$-IGDD of type $(11,1)^n$ for any $n\geq 3$ and $n\neq 6,8$.
\end{Lemma}

\proof For $n=3,4,5$, we construct a $(K_4-e)$-IGDD of type $(11,1)^n$ on
$I_{11n}$ with groups $\{nj+i: 0\leq j\leq 10\}$, $0\leq i\leq n-1$, and the hole set $I_{n}$. Let $\alpha_{n}=(0\ 1 \cdots n-1)(n\ n+1 \cdots 2n-1\ 6n\ 6n+1
 \cdots 7n-1)(2n\ 2n+1 \cdots 3n-1\ 7n\ 7n+1 \cdots 8n-1)\cdots(5n\ 5n+1 \cdots 6n-1\ 10n\ 10n+1 \cdots 11n-1)$ be a permutation on $I_{11n}$. Let $G_{n}$ be the group generated by $\alpha_{n}$. Only base blocks are listed below.

\begin{center}
\begin{tabular}{lllll}
$n=3$:&
$[0,5,4-7]$,&
$[0,10,8-11]$,&
$[13,14,0-3]$,&
$[16,17,0-3]$,\\&
$[3,  10,  5-  29]$,&
$[3,  11,  7-  31]$,&
$[8, 25,  3-  6]$,&
$[3,  32,  22-  28]$,\\&
$[7,  32,  6-  30]$,&
$[13,  29,  6-  21]$,&
$[11,  16,  9-  27]$,&
$[9,  31,  14-  29]$.\\
$n=4$:&
$[0,  6,  5-  9]$,&
$[0,  13,  7-  10]$,&
$[0, 17,  11-  14]$,&
$[0,  21,  15-  18]$,\\&
$[0,  42,  19-  23]$,&
$[4,  11,  6-  10]$,&
$[4,  14,  7-  13]$,&
$[4,  17,  15-  18]$,\\&
$[4,  21,  19-  22]$,&
$[4,  29,  23-  38]$,&
$[4,  37,  30-  39]$,&
$[33,  42,  4-  11]$,\\&
$[41,  43,  4-  8]$,&
$[8,  37,  10-  11]$,&
$[8,  34,  13-  21]$,&
$[8,  22,  15-  33]$,\\&
$[14,  19,  12-  40]$,&
$[18,  37,  12-  43]$.\\
$n=5$:&
$[0,  7, 6-  9]$,&
$[0,  11,  8-  12]$,&
$[0,  16,  13-  14]$,&
$[0,  18,  17-  21]$,\\&
$[0,  23,  19-  22]$,&
$[0,  26, 24-  27]$,&
$[28,  54,  0-  5]$,&
$[5,  12,  8-  14]$,\\&
$[5,  11,  9-  17]$,&
$[5,  21,  13-  19]$,&
$[5,  18,  16-  24]$,&
$[5,  26,  22-  23]$,\\&
$[5,  29,  27-  36]$,&
$[5,  44, 41-  51]$,&
$[42,  49,  5-  13]$,&
$[43,  52,  5-  19]$,\\&
$[46,  48,  5-  14]$,&
$[47,  53,  5-  14]$,&
$[10,  41,  13-  27]$,&
$[10,  26,  14-  44]$,\\&
$[17,  51,  10-  28]$,&
$[10,  47,  21-  54]$,&
$[28,  49, 10-  40]$,&
$[23,  46,  15-  54]$.\\
\end{tabular}
\end{center}
All other blocks are obtained by developing these base blocks under the action of $G_{n}$.

For $n\geq 7$ and $n\neq 8$, start from an $(n,\{3,4,5\},1)$-PBD (from Lemma \ref{lem:PBD(3,4,5)}), and then apply Construction \ref{con:PBD} with a $(K_4-e)$-IGDD of type $(11,1)^k$, $k\in\{3,4,5\}$, to obtain a $(K_4-e)$-IGDD of type $(11,1)^n$. \qed

\begin{Lemma}\label{lem:(11,6)^3}
There exists a $(K_4-e)$-IGDD of type $(11,6)^{n}$ for $n\equiv 1,3\pmod{6}$ and $n\geq 3$.
\end{Lemma}

\proof For $n=3$, we construct a $(K_4-e)$-IGDD of type $(11,6)^{3}$ on $I_{33}$ with groups $\{3j+i: 0\leq j\leq 10\}$, $0\leq i\leq 2$, and the hole set
$I_{18}$. Let $\alpha=\prod_{j=0}^{10}(3j\ 3j+1\ 3j+2)$ be a permutation on $I_{33}$. Let $G$ be the group generated by $\alpha$. Only base blocks are listed below:

\begin{center}
\begin{tabular}{lllll}
$[0,  20,  19-  25]$,&
$[22,  23,  0-  3]$,&
$[0,  28,  26-  29]$,&
$[31,  32,  0-  3]$,&
$[19,  29,  3-  6]$,\\
$[20,  28,  3-  9]$,&
$[25,  26,  3-  12]$,&
$[20,  22,  6-  30]$,&
$[23,  28,  6-  30]$,&
$[25,  32,  6-  9]$,\\
$[26,  31,  6-  9]$,&
$[19,  23,  9-  12]$,&
$[22,  29,  9-  12]$,&
$[20,  31,  12-  15]$,&
$[28,  32,  12-  15]$,\\
$[15,  26,  19-  22]$,&
$[15,  25,  23-  29]$.\\
\end{tabular}
\end{center}
All other blocks are obtained by developing these base blocks under the action of $G$.

For $n\equiv 1,3\pmod{6}$ and $n\geq 7$, start from an $(n,\{3\},1)$-PBD, which is also referred to as a Steiner triple system of order $n$ (cf. \cite{StinsonBook}), and then apply Construction \ref{con:PBD} with a $(K_4-e)$-IGDD of type $(11,6)^3$ to obtain a $(K_4-e)$-IGDD of type $(11,6)^n$. \qed

\section{Proof of Theorem \ref{thm:main}}

\subsection{The cases $g=2,3$ and $4$}

\begin{Lemma} \label{lem:2^4}
There exists a $(K_{4}-e)$-MGDP of type $2^4$ with leave $E_{4,j}$ for $2\leq j\leq 11$.
\end{Lemma}
\proof We give explicit constructions for the required MGDPs as follows:

\begin{center}
\begin{tabular}{ll}
$E_{4,2}=\{\{2,3\},\{3,6\},\{6,7\},\{7,2\}\}:$&
$[0,  1,  2-  3]$,
$[0,  5,  6-  7]$,
$[1,  4,  6-  7]$,
$[4,  5,  2-  3]$;\\

$E_{4,3}=\{\{1,3\},\{3,4\},\{4,7\},\{7,5\}\}:$&
$[0,  2,  1-  7]$,
$[0,  6,  3-  5]$,
$[1,  6,  4-  7]$,
$[2,  5,  3-  4]$;\\

$E_{4,4}=\{\{0,5\},\{0,6\},\{0,7\},\{2,5\}\}:$&
$[0,  3,  1-  2]$,
$[2,  4,  1-  7]$,
$[3,  4,  5-  6]$,
$[6,  7,  1-  5]$;\\

$E_{4,5}=\{\{0,1\},\{0,2\},\{0,3\},\{0,7\}\}:$&
$[1,  3,  2-  4]$,
$[2,  5,  4-  7]$,
$[5,  6,  0-  3]$,
$[6,  7,  1-  4]$;\\

$E_{4,6}=\{\{0,1\},\{0,2\},\{1,2\},\{4,5\}\}:$&
$[3,  4,  1-  2]$,
$[3,  6,  0-  5]$,
$[5,  7,  0-  2]$,
$[6,  7,  1-  4]$;\\

$E_{4,7}=\{\{0,1\},\{0,2\},\{0,3\},\{4,5\}\}:$&
$[0,  6,  5-  7]$,
$[1,  4,  2-  7]$,
$[2,  5,  3-  7]$,
$[3,  6,  1-  4]$;\\

$E_{4,8}=\{\{0,1\},\{1,2\},\{4,5\},\{5,6\}\}:$&
$[0,  3,  2-  5]$,
$[2,  7,  4-  5]$,
$[3,  4,  1-  6]$,
$[6,  7,  0-  1]$;\\

$E_{4,9}=\{\{2,7\},\{7,6\},\{6,3\},\{4,5\}\}:$&
$[0,  1,  2-  3]$,
$[0,  5,  6-  7]$,
$[1,  4,  6-  7]$,
$[2,  3,  4-  5]$;\\

$E_{4,10}=\{\{1,2\},\{2,3\},\{4,5\},\{6,7\}\}:$&
$[0,  1,  3-  6]$,
$[0,  5,  2-  7]$,
$[3,  6,  4-  5]$,
$[4,  7,  1-  2]$;\\

$E_{4,11}=\{\{0,1\},\{2,3\},\{4,5\},\{6,7\}\}:$&
$[0,  5,  2-  7]$,
$[1,  4,  3-  6]$,
$[2,  7,  1-  4]$,
$[3,  6,  0-  5]$;\\
\end{tabular}
\end{center}
where the vertex set is $I_8$ and the group set is $\{\{i,i+4\} :$ $0\leq i\leq 3\}$. \qed

\begin{Lemma} \label{lem:2^7}
There exists a $(K_{4}-e)$-MGDP of type $2^{7}$ with leave $E_{4}$.
\end{Lemma}

\proof We give an explicit constructions for a $(K_{4}-e)$-MGDP of type $2^{7}$ with leave $E_{4,1}$ as follows. Its vertex set is $I_{14}$, the group set is $\{\{i,i+7\} :$ $0\leq i\leq 6\}$, and $E_{4,1}$ consists of edges $\{9,10\},\{9,12\},\{10,12\}$ and $\{5,9\}$.

\begin{center}
\begin{tabular}{llllll}
$[0,8,6-13]$,&
$[1,9,6-13]$,&
$[2,10,6-13]$,&
$[3,11,6-13]$,&
$[4,12,6-13]$,&
$[5,7,6-13]$,\\
$[0,  2,  4-  5]$,&
$[0,  9,  3- 11]$,&
$[0,  1,  10- 12]$,&
$[1,  2,  7- 11]$,&
$[1,  5,  3-  4]$,&
$[2,  3,  8- 12]$,\\
$[4,  7,  3-  10]$,&
$[5,  10,  8- 11]$,&
$[8,  9,  4-  7]$,&
$[11, 12,  7-  8]$.\\
\end{tabular}
\end{center}
For other leaves, take a $(K_4-e)$-IGDD of type $2^{(7,4)}$ (from Lemma \ref{lem:1-factorization}), and then apply Construction \ref{con:filling by MGDP} with a $(K_{4}-e)$-MGDP of type $2^4$  with leave $E_{4,j}$ for $2\leq j\leq 11$ (from Lemma \ref{lem:2^4}). \qed

\begin{Lemma} \label{lem:2^8}
There exists a $(K_{4}-e)$-MGDP of type $2^8$ with leave $E_{2}$.
\end{Lemma}
\proof We construct a $(K_{4}-e)$-MGDP of type $2^{8}$ with leave $E_{2,1}$ as follows, where the vertex set is $I_{16}$, the group set is $\{\{i,i+8\} :$ $0\leq i\leq 7\}$, and $E_{2,1}$ consists of edges $\{0,1\}$ and $\{1,8\}$.

\begin{center}
\begin{tabular}{llllll}
$[0,  3,  2-  4]$,&
$[0,  6,  5-  7]$,&
$[0,  10,  9-  11]$,&
$[12,  13,  0-  1]$,&
$[14,  15,  0-  1]$,&
$[1,  4,  2-  5]$,\\
$[1,  6,  3-  10]$,&
$[7,  11,  1-  2]$,&
$[2,  8,  5-  6]$,&
$[2,  12,  9-  14]$,&
$[13,  15,  2-  3]$,&
$[3,  7,  5-  9]$,\\
$[9,  15,  5-  6]$,&
$[4,  10,  8-  15]$,&
$[10,  12,  3-  7]$,&
$[10,  14,  5-  13]$,&
$[4,  14,  9-  11]$,&
$[8,  14,  3-  7]$,\\
$[4,  13,  6-  7]$,&
$[11,  12,  5-  6]$,&
$[9,  13,  8- 11]$,&
$[8,  15,  11-  12]$.
\end{tabular}
\end{center}
Next start from a $(K_4-e)$-IGDD of type $2^{(8,3)}$ (from Lemma \ref{lem:IGDD g=2,3,4}), and then apply Construction \ref{con:filling by MGDP} with a $(K_{4}-e)$-MGDP of type $2^3$ with leave $E_{2,2}$ (from Example \ref{ex:2^3}). \qed

\begin{Lemma} \label{lem:2^9}
There exists a $(K_{4}-e)$-MGDP of type $2^9$ with leave $E_{4}$.
\end{Lemma}
\proof We construct a $(K_{4}-e)$-MGDP of type $2^{9}$ with leave $E_{4,1}$ as follows, where the vertex set is $I_{18}$, the group set is $\{\{i,i+9\} :$ $0\leq i\leq 8\}$, and $E_{4,1}$ consists of edges $\{0,1\}$, $\{0,2\}$, $\{1,2\}$ and $\{0,10\}$.

\begin{center}
{ \tabcolsep 0.04in
\begin{tabular}{llllll}
$[0,  4,  3-  5]$,&
$[0,  7,  6-  8]$,&
$[0,  12,  11-  13]$,&
$[14,  15,  0-  1]$,&
$[16,  17,  0-  1]$,&
$[1,  5,  3-  6]$,\\
$[1,  7,  4-  9]$,&
$[8,  12,  1-  2]$,&
$[11,  13,  1-  3]$,&
$[2,  6,  3-  4]$,&
$[2,  7,  5-  10]$,&
$[2,  13,  9-  16]$,\\
$[2,  17,  14-  15]$,&
$[3,  14,  7-  8]$,&
$[3,  10,  9-  17]$,&
$[15,  16,  3-  4]$,&
$[4,  9,  8-  17]$,&
$[10,  12,  4-  6]$,\\
$[11,  14,  4-  10]$,&
$[5,  15,  8-  12]$,&
$[9,  16,  5-  6]$,&
$[10,  13,  5-  15]$,&
$[11,  17,  5-  6]$,&
$[6,  13,  8-  14]$,\\
$[11,  15,  7-  9]$,&
$[7,  17,  12-  13]$,&
$[8,  16,  10-  11]$,&
$[12,  14,  9-  16]$.\\
\end{tabular}}
\end{center}
Next start from a $(K_4-e)$-IGDD of type $2^{(9,4)}$ (from Lemma \ref{lem:IGDD g=2,3,4}), and apply Construction \ref{con:filling by MGDP} with a $(K_{4}-e)$-MGDP of type $2^4$ with leave $E_{4,j}$ (from Lemma \ref{lem:2^4}), $2\leq j\leq 11$. \qed

\begin{Lemma} \label{lem:2^12}
There exists a $(K_{4}-e)$-MGDP of type $2^{12}$ with leave $E_4$.
\end{Lemma}

\proof We construct a $(K_{4}-e)$-MGDP of type $2^{12}$ with leave $E_{4,1}$ as follows, where the vertex set is $I_{24}$, the group set is $\{\{i,i+12\} :$ $0\leq i\leq 11\}$, and $E_{4,1}$ consists of edges $\{0,1\}$, $\{0,2\}$, $\{1,2\}$ and $\{0,3\}$.

\begin{center}
{\small \tabcolsep 0.05in
\begin{tabular}{llllll}
$[ 0, 4, 5-6 ]$,& $[ 0, 7, 8-9 ]$,& $[ 0, 10, 11-13 ]$,& $[ 0, 14, 15-16 ]$,& $[ 0, 17, 18-19 ]$,& $[ 20, 21, 0-1 ]$,\\
  $[ 22, 23, 0-1 ]$,& $[ 1, 3, 4-7 ]$,& $[ 1, 6, 5-17 ]$,& $[ 1, 8, 9-14 ]$,& $[ 1, 10, 12-15 ]$,& $[ 11, 16, 1-2 ]$,\\
  $[ 18, 19, 1-2 ]$,& $[ 2, 3, 5-10 ]$,& $[ 4, 8, 2-10 ]$,& $[ 7, 10, 5-17 ]$,& $[ 9, 10, 6-18 ]$,& $[ 10, 14, 19-20 ]$,\\
  $[ 10, 16, 21-23 ]$,& $[ 2, 9, 12-13 ]$,& $[ 2, 15, 6-17 ]$,& $[ 2, 7, 20-23 ]$,& $[ 21, 22, 2-3 ]$,& $[ 6, 7, 11-16 ]$,\\
  $[ 4, 7, 12-14 ]$,& $[ 7, 13, 15-21 ]$,& $[ 18, 22, 4-7 ]$,& $[ 4, 11, 9-21 ]$,& $[ 4, 13, 17-19]$,& $[ 4, 15, 20-23 ]$,\\
  $[ 5, 9, 14-15 ]$,& $[ 8, 15, 11-21 ]$,& $[ 12, 15, 16-18 ]$,& $[ 19, 22, 5-15 ]$,& $[ 5, 18, 8-21 ]$,& $[ 5, 13, 11-16 ]$,\\
  $[ 5, 12, 20-23 ]$,& $[ 3, 16, 8-18 ]$,& $[ 11, 14, 3-18 ]$,& $[ 13, 18, 20-23 ]$,& $[ 12, 19, 11-21 ]$,& $[ 11, 20, 17-22 ]$,\\
  $[ 16, 22, 9-17 ]$,& $[ 19, 20, 6-16 ]$,& $[ 8, 23, 17-19 ]$,& $[ 14, 17, 12-21 ]$,& $[ 6, 23, 14-21 ]$,& $[ 3, 9, 17-19 ]$,\\
  $[ 20, 23, 3-9 ]$,& $[ 13, 22, 8-14]$,& $[ 6, 12, 8-22 ]$,& $[ 3, 13, 6-12 ]$.\\
\end{tabular}}
\end{center}
For other leaves, start from a $(K_4-e)$-IGDD of type $2^{(12,4)}$ (from Lemma \ref{lem:g^(12,4)}) and apply Construction \ref{con:filling by MGDP} with a $(K_{4}-e)$-MGDP of type $2^4$ (from Lemma \ref{lem:2^4}). \qed

\begin{Lemma} \label{lem:2^13,2^14}
There exist a $(K_{4}-e)$-MGDP of type $2^{13}$ with leave $E_{2}$ and a $(K_{4}-e)$-MGDP of type $2^{14}$ with leave $E_{4}$.
\end{Lemma}

\proof By Lemmas \ref{lem:1-factorization} and \ref{lem:IGDD-2^(14,7)}, there exist a $(K_4-e)$-IGDD of type $2^{(13,8)}$ and a $(K_4-e)$-IGDD of type $2^{(14,7)}$. Apply Construction \ref{con:filling by MGDP} with a $(K_{4}-e)$-MGDP of type $2^{8}$ with leave $E_{2}$ (from Lemma \ref{lem:2^8}) and a $(K_{4}-e)$-MGDP of type $2^{7}$ with leave $E_{4}$ (from Lemma \ref{lem:2^7}), respectively. \qed

\begin{Lemma} \label{lem:3^4}
There exists a $(K_{4}-e)$-MGDP of type $3^4$ with leave $E_{4}$.
\end{Lemma}
\proof The required MGDP is constructed on $I_{9}$ with groups $\{i,i+4,i+8\}$, $0\leq i\leq 3$.

$E_{4,1}=\{\{0,1\},\{0,2\},\{1,2\},\{0,11\}\}:$

\begin{tabular}{llllll}
$[0,  3,  5-  6]$,&
$[0,  7,  9- 10]$,&
$[1,  3,  4- 10]$,&
$[2, 11,  4-  9]$,&
$[3,  8,  2-  9]$,&
$[4,  9,  6- 10]$,\\
$[5,  7,  2-  4]$,&
$[6, 11,  1-  5]$,&
$[7,  8,  1-  6]$,&
$[8, 10,  5- 11]$.\\
\end{tabular}

$E_{4,2}=\{\{0,1\},\{1,4\},\{4,5\},\{0,5\}\}:$

\begin{tabular}{llllll}
$[0,  3,  2-  9]$,&
$[6,  7,  0-  1]$,&
$[10,  11,  0-  9]$,&
$[2,  11,  1-  4]$,&
$[1,  10,  3-  8]$,&
$[5,  8,  2-  3]$,\\
$[7,  9,  2-  8]$,&
$[4,  6,  3-  9]$,&
$[7,  10,  4-  5]$,&
$[6,  11,  5-  8]$.\\
\end{tabular}

$E_{4,5}=\{\{0,1\},\{0,2\},\{0,3\},\{0,5\}\}:$

\begin{tabular}{llllll}
$[0,  6,  7- 9]$,&
$[1,  2,  3-  4]$,&
$[1,  8,  6-  7]$,&
$[2,  5,  7-  8]$,&
$[3,  4,  6-  9]$,&
$[3, 10,  5-  8]$,\\
$[5, 11,  4-  6]$,&
$[7, 10,  4-  9]$,&
$[9, 11,  2-  8]$,&
$[10, 11,  0-  1]$.\\
\end{tabular}

$E_{4,6}=\{\{1,4\},\{4,10\},\{1,10\},\{6,7\}\}:$

\begin{tabular}{llllll}
$[0,  1,  2-  3]$,&
$[0,  5,  6-  7]$,&
$[0, 10,  9- 11]$,&
$[2,  9,  3- 11]$,&
$[3,  5,  4- 10]$,&
$[4,  7,  2-  9]$,\\
$[5,  8,  2- 11]$,&
$[6,  8,  3-  9]$,&
$[6, 11,  1-  4]$,&
$[7,  8,  1- 10]$.\\
\end{tabular}

$E_{4,7}=\{\{0,1\},\{0,2\},\{0,3\},\{4,5\}\}:$

\begin{tabular}{llllll}
$[0,  5,  6-  7]$,&
$[0,  9, 10- 11]$,&
$[1,  2,  3-  7]$,&
$[1, 11,  8- 10]$,&
$[2,  5,  8- 11]$,&
$[3,  9,  6-  8]$,\\
$[3, 10,  4-  5]$,&
$[4,  6,  1- 11]$,&
$[4,  9,  2-  7]$,&
$[7,  8,  6- 10]$.\\
\end{tabular}

$E_{4,8}=\{\{2,7\},\{7,10\},\{4,9\},\{4,11\}\}:$

\begin{tabular}{llllll}
$[0,  1,  2-  3]$,&
$[0,  5,  6-  7]$,&
$[0,  9, 10- 11]$,&
$[1,  4,  6-  7]$,&
$[1, 10,  8- 11]$,&
$[2,  3,  4-  9]$,\\
$[2,  5,  8- 11]$,&
$[5, 10,  3-  4]$,&
$[6,  8,  3- 11]$,&
$[7,  9,  6-  8]$.\\
\end{tabular}

$E_{4,11}=\{\{0,1\},\{2,3\},\{4,5\},\{6,7\}\}:$

\begin{tabular}{llllll}
$[0,  2,  5-  7]$,&
$[0,  3,  9- 10]$,&
$[1,  4,  2- 11]$,&
$[1,  8,  7- 10]$,&
$[2,  8,  9- 11]$,&
$[3,  6,  1-  4]$,\\
$[4,  9,  7- 10]$,&
$[5,  8,  3-  6]$,&
$[5, 10,  7- 11]$,&
$[6, 11,  0-  9]$.\\
\end{tabular}

To obtain leave $E_{4,3}$, take the block $[2,3,4-9]$ from the $(K_{4}-e)$-MGDP of type $3^4$ with leave $E_{4,8}$, and replace it with a new block $[3,9,2-4]$. Then one can check that $E_{4,3}=\{\{2,4\},\{2,7\},\{7,10\},\{4,11\}\}$.

To obtain leave $E_{4,4}$, take the block $[3,10,4-5]$ from the $(K_{4}-e)$-MGDP of type $3^4$ with leave $E_{4,7}$, and replace it with a new block $[5,10,3-4]$. Then one can check that $E_{4,4}=\{\{0,1\},\{0,2\},\{0,3\},\{3,4\}\}$.

To obtain leave $E_{4,9}$, take the block $[2,  3,  4-  9]$ from the $(K_{4}-e)$-MGDP of type $3^4$ with leave $E_{4,8}$, and replace it with a new block $[4,  9,  2-  3]$. Then one can check that $E_{4,9}=\{\{2,3\},\{2,7\},\{7,10\},\{4,11\}\}$.

To obtain leave $E_{4,10}$, take the block $[2,  3,  4-  9]$ from the $(K_{4}-e)$-MGDP of type $3^4$ with leave $E_{4,8}$, and replace it with a new block $[2,  4,  3-  9]$. Then one can check that $E_{4,10}=\{\{3,9\},\{2,7\},\{7,10\},\{4,11\}\}$. \qed

\begin{Lemma} \label{lem 3^3,7}
There exist a $(K_{4}-e)$-MGDP of type $3^3$ with leave $E_{2,2}$ and a $(K_{4}-e)$-MGDP of type $3^{n}$ with leave listed in Table \ref{possible minimal leaves} for $n\in \{7,9,12,13,14\}$.
\end{Lemma}
\proof For $n\in \{3,7,9,13\}$, by Lemma \ref{lem:IGDD-(3,2)^n}, there exists a $(K_4-e)$-IGDD of type $(3,2)^{n}$. Start from this IGDD and then apply Construction \ref{con:fillling-2} with a $(K_{4}-e)$-MGDP of type $2^{n}$ (from Example \ref{ex:2^3}, Lemmas \ref{lem:2^7}, \ref{lem:2^9}-\ref{lem:2^13,2^14}).

For $n\in \{12,14\}$, start from a $(K_4-e)$-IGDD of type $3^{(n,4)}$ (from Lemma \ref{lem:g^(12,4)}) and apply Construction \ref{con:filling by MGDP} with a $(K_{4}-e)$-MGDP of type $3^4$ with leave $E_4$ (from Lemma \ref{lem:3^4}). \qed

\begin{Lemma} \label{(K4-e)-3^8}
There exists a $(K_{4}-e)$-MGDP of type $3^8$ with leave $E_{2}$.
\end{Lemma}

\proof We construct a $(K_{4}-e)$-MGDP of type $3^{8}$ with leave $E_{2,1}$ as follows, where the vertex set is $I_{24}$, the group set is $\{\{i,i+8,i+16\} :$ $0\leq i\leq 7\}$, and $E_{2,1}$ consists of edges $\{12,18\}$ and $\{12,22\}$.

\begin{center}
\begin{tabular}{lllll}
$[0,  2,  1-  3]$,&
$[0,  5,  4-  6]$,&
$[0,  9,  7-  10]$,&
$[0,  12,  11-  13]$,&
$[0,  15,  14-  17]$,\\
$[0,  19,  18-  20]$,&
$[0,  22,  21-  23]$,&
$[1,  4,  3-  6]$,&
$[1,  7,  5-  8]$,&
$[1,  11,  10-  13]$,\\
$[1,  14,  12-  16]$,&
$[1,  18,  15-  20]$,&
$[19,  22,  1-  2]$,&
$[21,  23,  1-  2]$,&
$[2,  7,  4-  6]$,\\
$[2,  8,  5-  9]$,&
$[2,  14,  11-  13]$,&
$[2,  15,  12-  16]$,&
$[17,  20,  2-  3]$,&
$[3,  9,  5-  6]$,\\
$[3,  10,  7-  8]$,&
$[3,  16,  12-  13]$,&
$[14,  21,  3-  4]$,&
$[15,  22,  3-  4]$,&
$[18,  23,  3-  4]$,\\
$[4,  11,  8-  9]$,&
$[4,  13,  10-  17]$,&
$[16,  19,  4-  5]$,&
$[5,  12,  10-  17]$,&
$[5,  15,  11-  20]$,\\
$[14,  23,  5-  8]$,&
$[18,  22,  5-  7]$,&
$[6,  13,  8-  15]$,&
$[6,  23,  10-  12]$,&
$[6,  18,  11-  16]$,\\
$[17,  19,  6-  14]$,&
$[20,  21,  6-  11]$,&
$[11,  17,  7-  23]$,&
$[7,  21,  12-  16]$,&
$[13,  19,  7-  23]$,\\
$[14,  20,  7-  10]$,&
$[12,  19,  8-  9]$,&
$[15,  21,  8-  9]$,&
$[17,  18,  8-  21]$,&
$[20,  22,  8-  13]$,\\
$[9,  18,  13-  14]$,&
$[16,  22,  9-  11]$,&
$[20,  23,  9-  16]$,&
$[10,  19,  15-  21]$,&
$[10,  17,  16-  22]$.\\
\end{tabular}
\end{center}
Next start from a $(K_4-e)$-IGDD of type $3^{(8,3)}$ (from Lemma \ref{lem:IGDD g=2,3,4}), and then apply Construction \ref{con:filling by MGDP} with a $(K_{4}-e)$-MGDP of type $3^3$ with leave $E_{2,2}$ (from Lemma \ref{lem 3^3,7}). \qed

\begin{Lemma} \label{lem:4^8}
There exists a $(K_{4}-e)$-MGDP of type $4^{8}$ with leave $E_{3}$.
\end{Lemma}
\proof We construct a $(K_{4}-e)$-MGDP of type $4^{8}$ with leave $E_{3,2}$ as follows, where the vertex set is $I_{32}$, the group set is $\{\{i,i+8,i+16,i+24\} :$ $0\leq i\leq 7\}$ and $E_{3,2}$ consists of edges $\{0,31\}$, $\{22,31\}$ and $\{30,31\}$. Its blocks are divided into two parts. The first part consists of blocks of a $(K_{4}-e)$-MGDP of type $2^{8}$, which is constructed on $I_{32}\setminus I_{16}$ with groups $\{i+16,i+24\}$,  $0\leq i\leq 7$. This MGDP exists by Lemma \ref{lem:2^8} and we can always make its leave contain edges $\{22,31\}$ and $\{30,31\}$. The second part consists of the following $67$ blocks:

\begin{center}
\begin{tabular}{llllll}
 $[ 1, 2, 0-31 ]$,& $[ 3, 6, 0-31 ]$,& $[ 4, 5, 0-31 ]$,& $[ 0, 9, 7-23 ]$,& $[ 0, 10, 17- 22 ]$,\\ $[ 0, 11, 18-30 ]$,&
  $[ 0, 13, 19-20 ]$,& $[ 0, 14, 21-28 ]$,& $[ 0, 12, 25-26 ]$,& $[ 0, 15, 27-29 ]$,\\ $[ 8, 9, 2-31 ]$,& $[ 10, 11, 1-31 ]$,&
  $[ 12, 31, 13-14 ]$,& $[ 3, 9, 4-29 ]$,& $[ 6, 9, 5-21 ]$,\\ $[ 9, 10, 16-27 ]$,& $[ 9, 12, 18-19 ]$,& $[ 9, 13, 22-30 ]$,&
  $[ 9, 11, 20-28 ]$,& $[ 14, 26, 1-9 ]$,\\ $[ 15, 24, 1-9 ]$,& $[ 1, 3, 5-21 ]$,& $[ 1, 6, 16-20 ]$,& $[ 1, 7, 18-28 ]$,&
  $[ 1, 4, 19-22 ]$,\\ $[ 1, 8, 23-29 ]$,& $[ 13, 27, 1-2 ]$,& $[ 12, 30, 1-2 ]$,& $[ 5, 12, 7-27 ]$,& $[ 6, 12, 8-29 ]$,\\
  $[ 4, 6, 2-27 ]$,& $[ 6, 13, 17-18 ]$,& $[ 6, 10, 19-28 ]$,& $[ 6, 11, 23-24 ]$,& $[ 7, 26, 3-6 ]$,\\ $[ 15, 25, 2-6 ]$,&
  $[ 3, 12, 16-23 ]$,& $[ 10, 24, 3-12 ]$,& $[ 12, 15, 21-22 ]$,& $[ 11, 17, 7-12 ]$,\\ $[ 2, 11, 21-29 ]$,& $[ 2, 3, 22-28 ]$,&
  $[ 7, 8, 20-22 ]$,& $[ 5, 11, 22-26 ]$,& $[ 14, 27, 7-8 ]$,\\ $[ 2, 7, 19-24 ]$,& $[ 5, 14, 19-24 ]$,& $[ 4, 13, 23-24 ]$,&
  $[ 2, 14, 20-23 ]$,& $[ 2, 5, 16-17 ]$,\\ $[ 5, 10, 20-23 ]$,& $[ 3, 15, 18-20 ]$,& $[ 5, 8, 18-25 ]$,& $[ 5, 15, 28-30 ]$,&
  $[ 8, 30, 3-10 ]$,\\ $[ 8, 13, 26-28 ]$,& $[ 8, 15, 11-19 ]$,& $[ 4, 8, 17-21 ]$,& $[ 14, 17, 3-15 ]$,& $[ 4, 15, 10-26 ]$,\\
  $[ 4, 7, 25-30 ]$,& $[ 7, 10, 21-29 ]$,& $[ 13, 25, 3-11 ]$,& $[ 10, 14, 13-25 ]$,& $[ 13, 16, 7-15 ]$,\\ $[ 4, 14, 18-29 ]$,&
  $[ 11, 16, 4-14 ]$.
\end{tabular}
\end{center}

To obtain a $(K_{4}-e)$-MGDP of type $4^8$ with leave $E_{3,1}$, begin with a $(K_{4}-e)$-IGDD of type $4^{(8,3)}$ (from Lemma \ref{lem:IGDD g=2,3,4}), and apply Construction \ref{con:filling by MGDP} with a $(K_{4}-e)$-MGDP of type $4^{3}$ (from Example \ref{eg:4^3}).

For other leaves, start from a $(K_4-e)$-IGDD of type $(4,1)^{8}$ (from Lemma \ref{lem:IGDD-(4,1)^5)}) and then apply Construction \ref{con:fillling-2} with a $(K_{4}-e)$-MGDP of type $1^8$ with leave $E_{3,3}$, $E_{3,4}$ and $E_{3,5}$ (from Theorem \ref{thm:1^n}). \qed

\begin{Lemma} \label{lem:(K4-e)-4^7}
There exists a $(K_{4}-e)$-MGDP of type $4^{n}$ with leave listed in Table \ref{possible minimal leaves} for $n\in \{7,9,12,13,14\}$.
\end{Lemma}

\proof For $n\in \{7,9,12,14\}$, start from a $(K_4-e)$-IGDD of type $4^{(n,4)}$ (from Lemmas \ref{lem:1-factorization}, \ref{lem:IGDD g=2,3,4} and \ref{lem:g^(12,4)}). Then apply Construction \ref{con:filling by MGDP} with a $(K_{4}-e)$-MGDP of type $4^4$ with leave $E_1$ (from Lemma \ref{eg:IGDD-(4,1)^4}).

For $n=13$, start from a $(K_4-e)$-IGDD of type $4^{(13,8)}$ (from Lemma \ref{lem:1-factorization}). Then apply Construction \ref{con:filling by MGDP} with a $(K_{4}-e)$-MGDP of type $4^8$ with leave $E_3$ (from Lemma \ref{lem:4^8}). \qed

\begin{Lemma} \label{lem:g=2,3,4}
Let $g\in\{2,3,4\}$. There exists a $(K_{4}-e)$-MGDP of type $g^n$ with leave listed in Table \ref{possible minimal leaves} for any $n\geq 4$, except when $(g,n)=(2,4)$.
\end{Lemma}

\proof For $n\equiv 0,1\pmod{5}$, the conclusion follows from Theorem \ref{thm:K4-e-GDD}. For $n\equiv 2,3,4\pmod{5}$ and $4\leq n\leq 14$, the conclusion follows from Example \ref{eg:IGDD-(4,1)^4} and Lemmas \ref{lem:2^4}-\ref{lem:(K4-e)-4^7}.

Let $h\in\{2,3,4\}$. By Lemma \ref{lem:IGDD g=2,3,4}, there exists a $(K_{4}-e)$-IGDD of type $g^{(5+h,h)}$ for $g\in\{2,3,4\}$. Then by Lemma \ref{lem:filling 1-type IGDD} there exists a $(K_{4}-e)$-IGDD of type $g^{(n,5+h)}$ for any $n\equiv h\pmod{5}$ and $n\geq 15+h$. Apply Construction \ref{con:filling by MGDP} with a $(K_{4}-e)$-MGDP of type $g^{5+h}$ to complete the proof. \qed

\subsection{The cases $g=6,7,8$ and $9$}

\begin{Lemma} \label{lem:6^n}
There exists a $(K_{4}-e)$-MGDP of type $6^{n}$ with leave listed in Table \ref{possible minimal leaves} for any $n\geq 3$.
\end{Lemma}

\proof For $n=3$, we construct a $(K_{4}-e)$-MGDP of type $6^{3}$ on $I_{18}$ with groups $\{3i+j:0\leq i\leq 5\}$, $0\leq j \leq 2$, as follows:

$E_{3,2}=\{\{9,7\},\{9,13\},\{9,16\}\}:$

\begin{tabular}{lllll}
$[0,  2,  1-  4]$,&
$[0,  7,  5-  8]$,&
$[0,  11,  10-  13]$,&
$[0,  16,  14-  17]$,&
$[1,  5,  3-  6]$,\\
$[1,  9,  8-  11]$,&
$[12,  14,  1-  4]$,&
$[15,  17,  1-  7]$,&
$[2,  13,  3-  6]$,&
$[7,  12,  2-  11]$,\\
$[9,  10,  2-  14]$,&
$[15,  16,  2-  11]$,&
$[4,  11,  3-  6]$,&
$[7,  14,  3-  6]$,&
$[8,  16,  3-  6]$,\\
$[10,  17,  3-  6]$,&
$[4,  9,  5-  17]$,&
$[8,  15,  4-  10]$,&
$[5,  12,  10-  16]$,&
$[13,  15,  5-  14]$,\\
$[12,  13,  8-  17]$.\\
\end{tabular}

$E_{3,3}=\{\{0,1\},\{1,2\},\{2,3\}\}:$

\begin{tabular}{lllll}
$[0,  4,  2-  5]$,&
$[0,  8,  7-  10]$,&
$[0,  13,  11-  14]$,&
$[16,  17,  0-  6]$,&
$[3,  11,  1-  7]$,\\
$[5,  9,  1-  10]$,&
$[6,  14,  1-  4]$,&
$[8,  12,  1-  4]$,&
$[15,  17,  1-  7]$,&
$[6,  10,  2-  11]$,\\
$[7,  9,  2-  14]$,&
$[12,  16,  2-  11]$,&
$[13,  15,  2-  5]$,&
$[3,  17,  4-  10]$,&
$[3,  16,  5-  14]$,\\
$[8,  13,  3-  6]$,&
$[4,  11,  9-  15]$,&
$[5,  7,  6-  12]$,&
$[8,  16,  9-  15]$,&
$[13,  17,  9-  12]$,\\
$[10,  14,  12-  15]$.\\
\end{tabular}

$E_{3,4}=\{\{9,13\},\{9,16\},\{10,14\}\}:$

\begin{tabular}{lllll}
$[0,  2,  1-  4]$,&
$[0,  7,  5-  8]$,&
$[0,  11,  10-  13]$,&
$[0,  16,  14-  17]$,&
$[1,  5,  3-  6]$,\\
$[1,  9,  8-  11]$,&
$[12,  14,  1-  4]$,&
$[15,  17,  1-  10]$,&
$[2,  10,  3-  9]$,&
$[6,  16,  2-  11]$,\\
$[7,  12,  2-  11]$,&
$[13,  15,  2-  5]$,&
$[4,  11,  3-  15]$,&
$[7,  17,  3-  6]$,&
$[8,  16,  3-  15]$,\\
$[13,  14,  3-  6]$,&
$[4,  9,  5-  17]$,&
$[6,  8,  4-  10]$,&
$[5,  12,  10-  16]$,&
$[7,  14,  9-  15]$,\\
$[12,  13,  8-  17]$.\\
\end{tabular}

$E_{3,5}=\{\{8,10\},\{9,11\},\{12,17\}\}:$

\begin{tabular}{lllll}
$[0,  2,  1-  4]$,&
$[0,  7,  5-  8]$,&
$[0,  11,  10-  13]$,&
$[0,  16,  14-  17]$,&
$[1,  5,  3-  6]$,\\
$[1,  9,  8-  14]$,&
$[11,  12,  1-  4]$,&
$[15,  17,  1-  10]$,&
$[2,  10,  3-  6]$,&
$[2,  15,  7-  13]$,\\
$[2,  16,  9-  12]$,&
$[4,  17,  3-  6]$,&
$[7,  11,  3-  6]$,&
$[8,  16,  3-  6]$,&
$[13,  14,  3-  6]$,\\
$[5,  9,  4-  10]$,&
$[4,  15,  8-  14]$,&
$[12,  13,  5-  8]$,&
$[15,  16,  5-  11]$,&
$[9,  17,  7-  13]$,\\
$[12,  14,  7-  10]$.\\
\end{tabular}

\noindent To obtain a $(K_{4}-e)$-MGDP of type $6^3$ with leave $E_{3,1}$, begin with a $(K_{4}-e)$-IGDD of type $(6,4)^{3}$ (from Lemma \ref{lem:IGDD-(3,2)^n}), and then apply Construction \ref{con:fillling-2} with a $(K_{4}-e)$-MGDP of type $4^{3}$ (from Example \ref{eg:4^3}).

For $n\geq 4$, by Lemma \ref{lem:IGDD-(3,2)^n}, there exists a $(K_{4}-e)$-IGDD of type $(6,4)^{n}$. Start from this IGDD, and apply Construction \ref{con:fillling-2} with a $(K_{4}-e)$-MGDP of type $4^{n}$ (from Lemma \ref{lem:g=2,3,4}). \qed

\begin{Lemma} \label{lem:7^n}
There exists a $(K_{4}-e)$-MGDP of type $7^{n}$ with leave listed in Table \ref{possible minimal leaves} for any $n\geq 3$.
\end{Lemma}

\proof For $n=3$, we construct a $(K_{4}-e)$-MGDP of type $7^{3}$ on $I_{21}$ with groups $\{3i+j:0\leq i\leq 6\}$, $0\leq j \leq 2$, as follows:
\begin{center}
{ \tabcolsep 0.03in
\begin{tabular}{llllll}
$[0,  2,  1-  4]$,&
$[0,  7,  5-  8]$,&
$[0,  11,  10-  13]$,&
$[0,  16,  14-  17]$,&
$[19,  20,  0-  3]$,&
$[1,  5,  3-  6]$,\\
$[1,  9,  8-  11]$,&
$[1,  14,  12-  15]$,&
$[1,  18,  17-  20]$,&
$[2,  10,  3-  12]$,&
$[2,  13,  6-  9]$,&
$[2,  18,  7-  19]$,\\
$[15,  16,  2-  8]$,&
$[3,  14,  4-  7]$,&
$[3,  13,  8-  17]$,&
$[11,  16,  3-  18]$,&
$[5,  18,  4-  13]$,&
$[4,  11,  6-  15]$,\\
$[4,  12,  8-  17]$,&
$[9,  20,  4-  7]$,&
$[5,  16,  9-  12]$,&
$[5,  15,  10-  19]$,&
$[7,  17,  6-  15]$,&
$[6,  19,  8-  14]$,\\
$[6,  20,  10-  16]$,&
$[11,  12,  7-  19]$,&
$[10,  18,  8-  14]$,&
$[9,  17,  10-  19]$,&
$[13,  20,  12-  15]$.\\
\end{tabular}
}
\end{center}
\noindent Its leave $E_{2,1}$ consists of two edges $\{9, 14\}$ and $\{13,14\}$. To obtain a $(K_{4}-e)$-MGDP of type $7^3$ with leave $E_{2,2}$, begin with a $(K_{4}-e)$-IGDD of type $(7,3)^{3}$ (from Lemma \ref{lem:(7,3)^n}), and then apply Construction \ref{con:fillling-2} with a $(K_{4}-e)$-MGDP of type $3^{3}$ (from Lemma \ref{lem 3^3,7}).

For $n=6$, the conclusion follows from Theorem \ref{thm:K4-e-GDD}. For $n\geq 4$ and $n\neq 6$, by Lemma \ref{lem:(7,3)^n}, there exists a $(K_4-e)$-IGDD of type $(7,3)^n$. Start from this IGDD, and apply Construction \ref{con:fillling-2} with a $(K_{4}-e)$-MGDP of type $3^{n}$ (from Lemma \ref{lem:g=2,3,4}). \qed

\begin{Lemma} \label{lem:8^n}
There exists a $(K_{4}-e)$-MGDP of type $8^{n}$ with leave listed in Table \ref{possible minimal leaves} for any $n\geq 3$.
\end{Lemma}

\proof For $n=3$, we construct a $(K_{4}-e)$-MGDP of type $8^{3}$ on $I_{24}$ with groups $\{3i+j:0\leq i\leq 7\}$, $0\leq j \leq 2$, as follows:

\begin{center}
{ \tabcolsep 0.03in
\begin{tabular}{llllll}
$[0,  2,  1-  4]$,&
$[0,  7,  5-  8]$,&
$[0,  11,  10-  13]$,&
$[0,  16,  14-  17]$,&
$[19,  20,  0- 3]$,&
$[22,  23,  0-  3]$,\\
$[1,  5,  3-  6]$,&
$[1,  9,  8-  11]$,&
$[1,  14,  12-  15]$,&
$[1,  18,  17-  20]$,&
$[21,  23,  1-  4]$,&
$[2,  7,  3-  6]$,\\
$[2,  10,  9-  15]$,&
$[12,  13,  2-  20]$,&
$[2,  18,  16-  22]$,&
$[19,  21,  2-  11]$,&
$[4,  14,  3- 6]$,&
$[8,  10,  3-  12]$,\\
$[11,  16,  3-  12]$,&
$[13,  17,  3-  6]$,&
$[5, 9,  4-  19]$,&
$[4,  18,  8-  11]$,&
$[12,  17,  4-  7]$,&
$[15,  20,  4-  22]$,\\
$[10,  18,  5-  23]$,&
$[5,  22,  12-  21]$,&
$[5,  15,  13-  16]$,&
$[6,  22,  8- 11]$,&
$[6,  20,  10-  16]$,&
$[19,  23,  6-  12]$,\\
$[7,  20,  9-  21]$,&
$[7,  15,  11-  23]$,&
$[14,  18,  7-  19]$,&
$[8,  21,  13-  16]$,&
$[15,  19,  8-  17]$,&
$[9,  23,  13-  16]$,\\
$[9,  22,  14-  17]$,&
$[10,  21,  14-  17]$.\\

\end{tabular}
}
\end{center}
\noindent Its leave $E_{2,1}$ consists of two edges $\{13, 14\}$ and $\{13,18\}$. To obtain a $(K_{4}-e)$-MGDP of type $8^3$ with leave $E_{2,2}$, begin with a $(K_{4}-e)$-IGDD of type $(8,3)^{3}$ (from Lemma \ref{lem:IGDD-(6,1)^n,(7,2)^n,(8,3)^n,(9,4)^n}), and then apply Construction \ref{con:fillling-2} with a $(K_{4}-e)$-MGDP of type $3^{3}$ (from Lemma \ref{lem 3^3,7}).

For $n=6$, the conclusion follows from Theorem \ref{thm:K4-e-GDD}. For $n\geq 4$ and $n\neq 6$, by Lemma \ref{lem:IGDD-(6,1)^n,(7,2)^n,(8,3)^n,(9,4)^n}, there exists a $(K_4-e)$-IGDD of type $(8,3)^n$. Start from this IGDD, and apply Construction \ref{con:fillling-2} with a $(K_{4}-e)$-MGDP of type $3^{n}$ (from Lemma \ref{lem:g=2,3,4}). \qed

\begin{Lemma} \label{lem:9^n}
There exists a $(K_{4}-e)$-MGDP of type $9^{n}$ with leave listed in Table \ref{possible minimal leaves} for any $n\geq 3$.
\end{Lemma}

\proof For $n=3$, take a $(K_4-e)$-IGDD of type $(9,6)^3$ (from Lemma \ref{lem:IGDD-(3,2)^n}) and then Construction \ref{con:fillling-2} with a $(K_{4}-e)$-MGDP of type $6^{3}$ (from Lemma \ref{lem:6^n}).

For $n=6$, the conclusion follows from Theorem \ref{thm:K4-e-GDD}. For $n\geq 4$ and $n\neq 6$, by Lemma \ref{lem:IGDD-(6,1)^n,(7,2)^n,(8,3)^n,(9,4)^n}, there exists a $(K_4-e)$-IGDD of type $(9,4)^n$. Start from this IGDD, and apply Construction \ref{con:fillling-2} with a $(K_{4}-e)$-MGDP of type $4^{n}$ (from Lemma \ref{lem:g=2,3,4}). \qed

\subsection{The cases $g=11,12,13$ and $14$}

\begin{Lemma} \label{lem:11^n}
There exists a $(K_{4}-e)$-MGDP of type $11^{n}$ with leave listed in Table \ref{possible minimal leaves} for any $n\geq 3$.
\end{Lemma}

\proof For $n=4$ or $n\geq 10$, start from a $(K_{4}-e)$-IGDD of type $(11,1)^{n}$ (from Lemma \ref{lem:(11,1)^n}), and then apply Construction \ref{con:fillling-2} with a $(K_{4}-e)$-MGDP of type $1^{n}$ (from Theorem \ref{thm:1^n}).

For $n=5,6$, the conclusion follows from Theorem \ref{thm:K4-e-GDD}. For $n=3,7,9$, start from a $(K_{4}-e)$-IGDD of type $(11,6)^{n}$ (from Lemma \ref{lem:(11,6)^3}), and then apply Construction \ref{con:fillling-2} with a $(K_{4}-e)$-MGDP of type $6^{n}$ (from Lemma \ref{lem:6^n}).

For $n=8$, start from a $(K_{4}-e)$-MGDP of type $11^{(8,3)}$ (from Lemma \ref{lem:11^(8,3)}), and then apply Construction \ref{con:filling by MGDP} with a $(K_{4}-e)$-MGDP of type $11^3$. \qed

\begin{Lemma} \label{lem:12^n,14^n}
Let $g\in\{12,14\}$. There exists a $(K_{4}-e)$-MGDP of type $g^{n}$ with leave listed in Table \ref{possible minimal leaves} for any $n\geq 3$.
\end{Lemma}

\proof For $n=6$, the conclusion follows from Theorem \ref{thm:K4-e-GDD}. For $n\geq 3$ and $n\neq 6$, set $(g,h)\in\{(6,4),(7,3)\}$. Apply Construction \ref{con:inflation-3} with a $(K_{4}-e)$-IGDD of type $(g,h)^{n}$ (from Lemmas \ref{lem:IGDD-(3,2)^n} and \ref{lem:(7,3)^n}) to obtain a $(K_{4}-e)$-IGDD of type $(2g,2h)^{n}$. Then apply Construction \ref{con:fillling-2} with a $(K_{4}-e)$-MGDP of type $(2h)^{n}$ (from Lemmas \ref{lem:8^n} and \ref{lem:6^n}). \qed

\begin{Lemma} \label{lem:13^n}
There exists a $(K_{4}-e)$-MGDP of type $13^{n}$ with leave listed in Table \ref{possible minimal leaves} for any $n\geq 3$.
\end{Lemma}

\proof For $n=6$, the conclusion follows from Theorem \ref{thm:K4-e-GDD}. For $n\geq 3$ and $n\neq 6$, start from a $(K_{4}-e)$-HGDD of type $(n,3^2 6^1)$ (from Lemma \ref{lem:HGDD-(n,5^u)}(2)), and apply Construction \ref{con:filling HGDD} with a $(K_{4}-e)$-IGDD of type $(4,1)^{n}$ (from Lemma \ref{lem:IGDD-(4,1)^5)}) to obtain a $(K_{4}-e)$-IGDD of type $(13,7)^{n}$. Finally apply Construction \ref{con:fillling-2} with a $(K_{4}-e)$-MGDP of type $7^{n}$ (from Lemma \ref{lem:7^n}). \qed

\subsection{General $g$}

\begin{Lemma} \label{lem:general g}
There exists a $(K_{4}-e)$-MGDP of type $g^{n}$ with leave listed in Table \ref{possible minimal leaves} for any $g\equiv 1,2,3,4\pmod{5}$, $g\geq 16$ and $n\geq 3$.
\end{Lemma}

\proof For $n=6$, the conclusion follows from Theorem \ref{thm:K4-e-GDD}. For $n\geq 3$ and $n\neq 6$, by Lemma \ref{lem:IGDD-(6,1)^n,(7,2)^n,(8,3)^n,(9,4)^n}, there exists a $(K_4-e)$-IGDD of type $(5+r,r)^n$, $r\in\{1,2,3,4\}$. Then apply Lemma \ref{lem:filling 2-type IGDD} to obtain a $(K_{4}-e)$-IGDD of type $(g,5+r)^{n}$ for any $g\equiv r\pmod{5}$ and $g\geq 15+r$. Finally apply Construction \ref{con:fillling-2} with a $(K_{4}-e)$-MGDP of type $(5+r)^n$ (from Lemmas \ref{lem:6^n}-\ref{lem:9^n}). \qed

\noindent \textbf{Proof of Theorem \ref{thm:main}} When $g=1$, the conclusion follows from Theorem \ref{thm:1^n}. When $g\equiv 0\pmod{5}$ or $n\equiv 0,1\pmod{5}$, the conclusion follows from Theorem \ref{thm:K4-e-GDD}. When $g\in\{2,3,4,6,7,8,9,11,12,13,14\}$, the conclusion follows from Examples \ref{ex:2^3} and \ref{eg:4^3}, Lemmas \ref{lem:2^3_no}-\ref{lem:4^3_no}, \ref{lem:2^4}, \ref{lem 3^3,7} and Lemmas \ref{lem:g=2,3,4}-\ref{lem:13^n}. When $g\equiv 1,2,3,4\pmod{5}$ and $g\geq 16$, the conclusion follows from Lemma \ref{lem:general g}. \qed

\section{Remarks}

Maximum packings of graphs have been and continue to be popular topics of research. Maximum packings of $K_n$ have attracted the most attention (see \cite{abb,abz,bs,b,br,bh,e,hor} for example). Instead, very little has been done on maximum packings of $K_{n(g)}$. Especially less attention has been paid on the issue of the existence of maximum packings with ANY possible minimum leave.

Hu, Chang and Feng \cite{hcf} investigated all possible leaves of maximum group divisible packings with $K_3$ or $K_3+e$. This paper continues the study in \cite{hcf} to examine all possible minimum leaves for group divisible $(K_4-e)$-packings. The next challenging problem is to examine all possible minimum leaves for maximum group divisible packings with $K_4$. As pointed out in \cite{hcf}, this problem is closely related with many  combinatorial configurations, such as balanced sampling plans excluding contiguous units \cite{cl} and matching divisible design \cite{dfl}. These relationships among kinds of designs motivate the study of this problem.

Finally we remark that recently Dyer, Haghshenas and Shalaby \cite{dhs} examined all possible minimum leaves for maximum packings of $K_n$ with stars having up to five edges. A natural and straightforward task is to generalize their result to the case of $K_{n(g)}$.

\end{document}